\documentclass[11pt,a4paper]{article}

\usepackage[utf8]{inputenc}
\usepackage[T1]{fontenc}
\usepackage[margin=1in]{geometry}
\usepackage{amsmath}
\usepackage{amssymb}
\usepackage{stmaryrd}
\usepackage{booktabs}
\usepackage{array}
\usepackage{tabularx}
\usepackage{authblk}
\usepackage{xcolor}
\usepackage{graphicx}
\usepackage[export]{adjustbox}
\usepackage{longtable}
\usepackage{calc}
\usepackage{newunicodechar}
\usepackage{listings}
\usepackage{soul}
\usepackage[colorlinks=true,linkcolor=blue!40!black,citecolor=blue!40!black,urlcolor=blue!40!black]{hyperref}

\let\origtextunderscore\_
\renewcommand{\_}{\origtextunderscore\allowbreak}

\lstdefinestyle{leanbox}{
  basicstyle=\footnotesize\ttfamily,
  columns=fullflexible,
  keepspaces=true,
  breaklines=true,
  breakatwhitespace=false,
  showstringspaces=false,
  upquote=true,
  extendedchars=true,
  inputencoding=utf8,
  literate=
    {ä}{{\"{a}}}1
    {é}{{\'{e}}}1
    {∑}{{\ensuremath{\sum}}}1
    {‹}{{`}}1
    {›}{{'}}1
    {∈}{{\ensuremath{\in}}}1
    {·}{{\ensuremath{\cdot}}}1
    {ℝ}{{\ensuremath{\mathbb{R}}}}1
    {ℕ}{{\ensuremath{\mathbb{N}}}}1
    {∼}{{\ensuremath{\sim}}}1
    {≃}{{\ensuremath{\simeq}}}1
    {§}{{\S}}1
    {¬}{{\ensuremath{\neg}}}1
    {·}{{\ensuremath{\cdot}}}1
    {×}{{\ensuremath{\times}}}1
    {Γ}{{\ensuremath{\Gamma}}}1
    {Δ}{{\ensuremath{\Delta}}}1
    {Η}{{\ensuremath{H}}}1
    {Θ}{{\ensuremath{\Theta}}}1
    {Κ}{{\ensuremath{K}}}1
    {Σ}{{\ensuremath{\Sigma}}}1
    {α}{{\ensuremath{\alpha}}}1
    {β}{{\ensuremath{\beta}}}1
    {λ}{{\ensuremath{\lambda}}}1
    {σ}{{\ensuremath{\sigma}}}1
    {τ}{{\ensuremath{\tau}}}1
    {φ}{{\ensuremath{\varphi}}}1
    {χ}{{\ensuremath{\chi}}}1
    {ψ}{{\ensuremath{\psi}}}1
    {ᵢ}{{\textsubscript{i}}}1
    {ⱼ}{{\textsubscript{j}}}1
    {—}{{\textemdash}}1
    {–}{{\textendash}}1
    {…}{{\ldots}}1
    {′}{{\ensuremath{^{\prime}}}}1
    {⁺}{{\textsuperscript{+}}}1
    {⁻}{{\textsuperscript{\ensuremath{-}}}}1
    {ⁱ}{{\textsuperscript{i}}}1
    {¹}{{\textsuperscript{1}}}1
    {₀}{{\textsubscript{0}}}1
    {₁}{{\textsubscript{1}}}1
    {₂}{{\textsubscript{2}}}1
    {₃}{{\textsubscript{3}}}1
    {₊}{{\textsubscript{+}}}1
    {ₙ}{{\textsubscript{n}}}1
    {ₚ}{{\textsubscript{p}}}1
    {ₖ}{{\textsubscript{k}}}1
    {ₘ}{{\textsubscript{m}}}1
    {ₜ}{{\textsubscript{t}}}1
    {ι}{{\ensuremath{\iota}}}1
    {ω}{{\ensuremath{\omega}}}1
    {ḡ}{{\={g}}}1
    {ē}{{\={e}}}1
    {ń}{{\'{n}}}1
    {ˢ}{{\textsuperscript{s}}}1
    {⦃}{{\textbraceleft\textbraceleft}}2
    {⦄}{{\textbraceright\textbraceright}}2
    {𝒵}{{\ensuremath{\mathcal{Z}}}}1
    {𝓕}{{\ensuremath{\mathcal{F}}}}1
    {ℕ}{{\ensuremath{\mathbb{N}}}}1
    {←}{{\ensuremath{\leftarrow}}}1
    {↑}{{\ensuremath{\uparrow}}}1
    {→}{{\ensuremath{\rightarrow}}}1
    {↓}{{\ensuremath{\downarrow}}}1
    {↔}{{\ensuremath{\leftrightarrow}}}1
    {↥}{{\ensuremath{\Uparrow}}}1
    {↦}{{\ensuremath{\mapsto}}}1
    {↪}{{\ensuremath{\hookrightarrow}}}1
    {↟}{{\ensuremath{\Uparrow}}}1
    {⇄}{{\ensuremath{\rightleftarrows}}}1
    {⇑}{{\ensuremath{\Uparrow}}}1
    {⇒}{{\ensuremath{\Rightarrow}}}1
    {⨆}{{\ensuremath{\bigcup}}}1
    {𝕆}{{\ensuremath{\mathbb{O}}}}1
    {∀}{{\ensuremath{\forall}}}1
    {∃}{{\ensuremath{\exists}}}1
    {∅}{{\ensuremath{\emptyset}}}1
    {∈}{{\ensuremath{\in}}}1
    {∉}{{\ensuremath{\notin}}}1
    {∘}{{\ensuremath{\circ}}}1
    {∋}{{\ensuremath{\ni}}}1
    {∩}{{\ensuremath{\cap}}}1
    {∏}{{\ensuremath{\prod}}}1
    {∞}{{\ensuremath{\infty}}}1
    {⋂}{{\ensuremath{\bigcap}}}1
    {⋃}{{\ensuremath{\bigcup}}}1
    {∧}{{\ensuremath{\wedge}}}1
    {∨}{{\ensuremath{\vee}}}1
    {∪}{{\ensuremath{\cup}}}1
    {∼}{{\ensuremath{\sim}}}1
    {≠}{{\ensuremath{\neq}}}1
    {≤}{{\ensuremath{\leq}}}1
    {≥}{{\ensuremath{\geq}}}1
    {≪}{{\ensuremath{\ll}}}1
    {≃}{{\ensuremath{\simeq}}}1
    {≅}{{\ensuremath{\cong}}}1
    {⊆}{{\ensuremath{\subseteq}}}1
    {⊑}{{\ensuremath{\sqsubseteq}}}1
    {⊔}{{\ensuremath{\sqcup}}}1
    {⊓}{{\ensuremath{\sqcap}}}1
    {⊢}{{\ensuremath{\vdash}}}1
    {⊣}{{\ensuremath{\dashv}}}1
    {⊤}{{\ensuremath{\top}}}1
    {⊥}{{\ensuremath{\bot}}}1
    {⊨}{{\ensuremath{\vDash}}}1
    {⊬}{{\ensuremath{\nvdash}}}1
    {⊭}{{\ensuremath{\nvDash}}}1
    {⋀}{{\ensuremath{\bigwedge}}}1
    {⋁}{{\ensuremath{\bigvee}}}1
    {─}{{\ensuremath{-}}}1
    {│}{{\ensuremath{|}}}1
    {└}{{+}}1
    {├}{{+}}1
    {□}{{\ensuremath{\square}}}1
    {▸}{{\ensuremath{\blacktriangleright}}}1
    {◇}{{\ensuremath{\lozenge}}}1
    {⟨}{{\ensuremath{\langle}}}1
    {⟩}{{\ensuremath{\rangle}}}1
    {⟶}{{\ensuremath{\longrightarrow}}}2
    {⟷}{{\ensuremath{\longleftrightarrow}}}2
    {⟹}{{\ensuremath{\Longrightarrow}}}2
}
\newunicodechar{α}{\ensuremath{\alpha}}
\newunicodechar{β}{\ensuremath{\beta}}
\newunicodechar{σ}{\ensuremath{\sigma}}
\newunicodechar{φ}{\ensuremath{\varphi}}
\newunicodechar{χ}{\ensuremath{\chi}}
\newunicodechar{ψ}{\ensuremath{\psi}}
\newunicodechar{Γ}{\ensuremath{\Gamma}}
\newunicodechar{Δ}{\ensuremath{\Delta}}
\newunicodechar{ℕ}{\ensuremath{\mathbb{N}}}
\newunicodechar{∀}{\ensuremath{\forall}}
\newunicodechar{∃}{\ensuremath{\exists}}
\newunicodechar{¬}{\ensuremath{\neg}}
\newunicodechar{∧}{\ensuremath{\wedge}}
\newunicodechar{□}{\ensuremath{\square}}
\newunicodechar{◇}{\ensuremath{\lozenge}}
\newunicodechar{∼}{\ensuremath{\sim}}
\newunicodechar{→}{\ensuremath{\rightarrow}}
\newunicodechar{←}{\ensuremath{\leftarrow}}
\newunicodechar{↔}{\ensuremath{\leftrightarrow}}
\newunicodechar{↦}{\ensuremath{\mapsto}}
\newunicodechar{⇒}{\ensuremath{\Rightarrow}}
\newunicodechar{⟶}{\ensuremath{\longrightarrow}}
\newunicodechar{⟷}{\ensuremath{\longleftrightarrow}}
\newunicodechar{⟹}{\ensuremath{\Longrightarrow}}
\newunicodechar{⟺}{\ensuremath{\Longleftrightarrow}}
\newunicodechar{∈}{\ensuremath{\in}}
\newunicodechar{⊢}{\ensuremath{\vdash}}
\newunicodechar{⊨}{\ensuremath{\vDash}}
\newunicodechar{≠}{\ensuremath{\neq}}
\newunicodechar{≈}{\ensuremath{\approx}}
\newunicodechar{≅}{\ensuremath{\cong}}
\newunicodechar{≪}{\ensuremath{\ll}}
\newunicodechar{⊆}{\ensuremath{\subseteq}}
\newunicodechar{⊇}{\ensuremath{\supseteq}}
\newunicodechar{⊋}{\ensuremath{\supsetneq}}
\newunicodechar{⊑}{\ensuremath{\sqsubseteq}}
\newunicodechar{⊔}{\ensuremath{\sqcup}}
\newunicodechar{⊓}{\ensuremath{\sqcap}}
\newunicodechar{⊤}{\ensuremath{\top}}
\newunicodechar{⊥}{\ensuremath{\bot}}
\newunicodechar{⊣}{\ensuremath{\dashv}}
\newunicodechar{↪}{\ensuremath{\hookrightarrow}}
\newunicodechar{↟}{\ensuremath{\Uparrow}}
\newunicodechar{⇄}{\ensuremath{\rightleftarrows}}
\newunicodechar{⇑}{\ensuremath{\Uparrow}}
\newunicodechar{⨆}{\ensuremath{\bigcup}}
\newunicodechar{𝕆}{\ensuremath{\mathbb{O}}}
\newunicodechar{∪}{\ensuremath{\cup}}
\newunicodechar{⊕}{\ensuremath{\oplus}}
\newunicodechar{×}{\ensuremath{\times}}
\newunicodechar{⟨}{\ensuremath{\langle}}
\newunicodechar{⟩}{\ensuremath{\rangle}}
\newunicodechar{₀}{\textsubscript{0}}
\newunicodechar{₁}{\textsubscript{1}}
\newunicodechar{ₙ}{\textsubscript{n}}
\newunicodechar{⁺}{\textsuperscript{+}}
\newunicodechar{′}{\ensuremath{^{\prime}}}
\newunicodechar{¹}{\textsuperscript{1}}
\newunicodechar{ⁱ}{\textsuperscript{i}}
\newunicodechar{⁻}{\textsuperscript{-}}
\newunicodechar{ⁿ}{\textsuperscript{n}}
\newunicodechar{₂}{\textsubscript{2}}
\newunicodechar{₃}{\textsubscript{3}}
\newunicodechar{ₖ}{\textsubscript{k}}
\newunicodechar{ₘ}{\textsubscript{m}}
\newunicodechar{ₚ}{\textsubscript{p}}
\newunicodechar{ₜ}{\textsubscript{t}}
\newunicodechar{ᵢ}{\textsubscript{i}}
\newunicodechar{ⱼ}{\textsubscript{j}}
\newunicodechar{₊}{\textsubscript{+}}
\newunicodechar{Η}{\ensuremath{\mathrm{H}}}
\newunicodechar{Θ}{\ensuremath{\Theta}}
\newunicodechar{Κ}{\ensuremath{\mathrm{K}}}
\newunicodechar{Σ}{\ensuremath{\Sigma}}
\newunicodechar{ι}{\ensuremath{\iota}}
\newunicodechar{λ}{\ensuremath{\lambda}}
\newunicodechar{τ}{\ensuremath{\tau}}
\newunicodechar{ω}{\ensuremath{\omega}}
\newunicodechar{≤}{\ensuremath{\leq}}
\newunicodechar{≥}{\ensuremath{\geq}}
\newunicodechar{≃}{\ensuremath{\simeq}}
\newunicodechar{∉}{\ensuremath{\notin}}
\newunicodechar{∋}{\ensuremath{\ni}}
\newunicodechar{∅}{\ensuremath{\emptyset}}
\newunicodechar{∘}{\ensuremath{\circ}}
\newunicodechar{∞}{\ensuremath{\infty}}
\newunicodechar{∨}{\ensuremath{\vee}}
\newunicodechar{∩}{\ensuremath{\cap}}
\newunicodechar{∏}{\ensuremath{\prod}}
\newunicodechar{⋀}{\ensuremath{\bigwedge}}
\newunicodechar{⋁}{\ensuremath{\bigvee}}
\newunicodechar{⋂}{\ensuremath{\bigcap}}
\newunicodechar{⊬}{\ensuremath{\nvdash}}
\newunicodechar{⊭}{\ensuremath{\nvDash}}
\newunicodechar{↑}{\ensuremath{\uparrow}}
\newunicodechar{↓}{\ensuremath{\downarrow}}
\newunicodechar{↥}{\ensuremath{\Uparrow}}
\newunicodechar{⟸}{\ensuremath{\Longleftarrow}}
\newunicodechar{⨅}{\ensuremath{\bigcap}}
\newunicodechar{◁}{\ensuremath{\triangleleft}}
\newunicodechar{▸}{\ensuremath{\blacktriangleright}}
\newunicodechar{𝓕}{\ensuremath{\mathcal{F}}}
\newunicodechar{─}{\-}
\newunicodechar{└}{\textasciigrave{}|}
\newunicodechar{├}{\textbar{}}
\newunicodechar{ē}{\={e}}
\newunicodechar{ń}{\'{n}}
\newunicodechar{ḡ}{\={g}}
\newunicodechar{·}{\textperiodcentered{}}
\newunicodechar{—}{\textemdash{}}
\newunicodechar{–}{\textendash{}}
\newunicodechar{…}{\ldots}
\newunicodechar{§}{\S}
\newunicodechar{│}{\ensuremath{\mid}}
\newunicodechar{✓}{\ensuremath{\checkmark}}
\newunicodechar{æ}{\ae}
\newunicodechar{ö}{\"o}

\newunicodechar{ä}{\"{a}}
\newunicodechar{é}{\'{e}}
\newunicodechar{∑}{\ensuremath{\sum}}
\newunicodechar{‹}{`}
\newunicodechar{›}{'}
\newunicodechar{∈}{\ensuremath{\in}}
\newunicodechar{·}{\ensuremath{\cdot}}
\newunicodechar{ℝ}{\ensuremath{\mathbb{R}}}
\newunicodechar{ℕ}{\ensuremath{\mathbb{N}}}
\newunicodechar{∼}{\ensuremath{\sim}}
\newunicodechar{≃}{\ensuremath{\simeq}}
\newunicodechar{⦃}{\textbraceleft\textbraceleft}
\newunicodechar{⦄}{\textbraceright\textbraceright}

\title{\textbf{On the Number of Distinct Topological Bases of a Finite Set of Size \(N\)}}

\author[1]{\textbf{Lars Warren Ericson}}
\affil[1]{ORCID: 0000-0001-8299-9361}
\affil[1]{Catskills Research Company}
\affil[1]{\texttt{lars.ericson@catskillsresearch.com}}

\date{\today}

\begin{document}
\maketitle

\begin{center}
  \small
  \textbf{Github:} \url{https://github.com/catskillsresearch/cardb} \\
  \textbf{Palomar Registration:} \url{https://palomar-registry.org/entry?id=PALOMAR-2026-08-20-000003&version=1}
\end{center}

\begin{abstract}
For a finite set \(S\) with \(\lvert S\rvert = N\), the number of families \(\mathcal{B} \subseteq \mathcal{P}(S)\) that are topological bases is
\(\#(N) = \sum_{\mathcal{T} \in \operatorname{Top}(S)} 2^{\lvert\mathcal{T}\rvert - \lvert\mathcal{M}_{\mathcal{T}}\rvert}\),
where \(\mathcal{M}_{\mathcal{T}}\) is the canonical minimal basis of minimal open neighborhoods. The identity is proved in Lean 4 / Mathlib (\texttt{CARDB.lean}): bases generating \(\mathcal{T}\) are exactly the sets with \(\mathcal{M}_{\mathcal{T}} \subseteq \mathcal{B} \subseteq \mathcal{T}\). The small-\(N\) table and the discrete-dominance sandwich are proved in \texttt{CARDB/\allowbreak{}SmallN.lean} and \texttt{CARDB/\allowbreak{}Asymptotics.lean}. Companion artifacts: those three Lean files (standalone Lake package in this directory).
\end{abstract}

\tableofcontents
\newpage

\hypertarget{introduction-and-definitions}{%
\section{Introduction and Definitions}\label{introduction-and-definitions}}

Let \(S\) be a set. In point-set topology, a family \(\mathcal{B} \subseteq \mathcal{P}(S)\) is defined as a \textbf{basis for a topology on \(S\)} if it satisfies two axioms:

\begin{enumerate}
\def\labelenumi{\arabic{enumi}.}
\item
  \textbf{Covering Axiom:}
  \[\bigcup_{B \in \mathcal{B}} B = S \quad \Longleftrightarrow \quad \forall x \in S,\; \exists B \in \mathcal{B} \text{ such that } x \in B\]
\item
  \textbf{Intersection Axiom:}
  \[\forall B_1, B_2 \in \mathcal{B},\; \forall x \in B_1 \cap B_2,\; \exists B_3 \in \mathcal{B} \text{ such that } x \in B_3 \subseteq B_1 \cap B_2\]
\end{enumerate}

The topology generated by \(\mathcal{B}\), written \(\operatorname{gen}(\mathcal{B})\) (Mathlib: \texttt{generateFrom}), is the collection of all arbitrary unions of subfamilies of \(\mathcal{B}\):
\[\operatorname{gen}(\mathcal{B}) = \left\{ \bigcup \mathcal{C} \;\middle|\; \mathcal{C} \subseteq \mathcal{B} \right\}\]
We reserve \(\mathcal{T}\) for a topology on \(S\), never for this generation map.

Companion files: \texttt{CARDB.lean}, \texttt{CARDB/\allowbreak{}SmallN.lean}, \texttt{CARDB/\allowbreak{}Asymptotics.lean}.

\begin{lstlisting}
import Mathlib.Topology.AlexandrovDiscrete
import Mathlib.Data.Fintype.BigOperators
import Mathlib.Data.Fintype.Powerset
import Mathlib.Data.Set.Card
import Mathlib.Data.Set.Finite.Basic
import Mathlib.Logic.Equiv.Sum

open Set TopologicalSpace
open scoped BigOperators

variable {α : Type*} [Fintype α] [DecidableEq α]

/-- A collection of subsets `B` is a topological basis if it covers `α`
    and satisfies the local intersection property. -/
def IsBasis (B : Set (Set α)) : Prop :=
  (⋃₀ B = univ) ∧
  ∀ ⦃U V⦄, U ∈ B → V ∈ B → ∀ x ∈ U ∩ V, ∃ W ∈ B, x ∈ W ∧ W ⊆ U ∩ V
\end{lstlisting}

\hypertarget{set-theoretic-status}{%
\subsection{Set-Theoretic Status}\label{set-theoretic-status}}

The definition of a basis and the topology generated by it is strictly constructible within standard Zermelo--Fraenkel set theory (\(\mathsf{ZF}\)) without the Axiom of Choice:
- \(\mathcal{P}(S)\) exists by the \textbf{Power Set Axiom}.
- \(\mathcal{B}\) exists as a subset of \(\mathcal{P}(S)\) by \textbf{Separation}.
- For each \(\mathcal{C} \in \mathcal{P}(\mathcal{B})\), \(\bigcup \mathcal{C}\) exists by the \textbf{Union Axiom}.
- The collection \(\operatorname{gen}(\mathcal{B})\) is the image of \(\mathcal{C} \mapsto \bigcup \mathcal{C}\) under the \textbf{Axiom of Replacement}.

\begin{center}\rule{0.5\linewidth}{0.5pt}\end{center}

\hypertarget{the-enumeration-problem}{%
\section{The Enumeration Problem}\label{the-enumeration-problem}}

Let \(S\) be a finite set with \(|S| = N\). We wish to determine the exact quantity:
\[\#(S) = \big| \left\{ \mathcal{B} \subseteq \mathcal{P}(S) \;\middle|\; \mathcal{B} \text{ is a topological basis on } S \right\} \big|\]

The objects being counted are the subtype of families satisfying \texttt{IsBasis}. Because \(\alpha\) is finite, there are only finitely many such families, and only finitely many topologies: a topology is recovered from its set of opens, and that assignment is injective.

\begin{lstlisting}
def ValidBasis (α : Type*) [Fintype α] := { B : Set (Set α) // IsBasis B }

/-- Open sets of `t`, as a family of subsets. This map is injective, so there
    are only finitely many topologies on a finite type. -/
def opensOf (t : TopologicalSpace α) : Set (Set α) := {U | t.IsOpen U}

omit [Fintype α] [DecidableEq α] in
theorem opensOf_injective :
    Function.Injective (opensOf : TopologicalSpace α → Set (Set α)) :=
  leftInverse_generateFrom.injective

noncomputable instance : DecidableEq (TopologicalSpace α) :=
  opensOf_injective.decidableEq

noncomputable instance : Fintype (TopologicalSpace α) :=
  Fintype.ofInjective opensOf opensOf_injective

noncomputable instance : DecidablePred (IsBasis : Set (Set α) → Prop) :=
  fun _ => Classical.dec _

noncomputable instance : Fintype (ValidBasis α) :=
  Subtype.fintype _
\end{lstlisting}

At first inspection, this appears to be a local hypergraph enumeration problem. However, the basis axioms enforce global closure properties that reduce the count to the classification of finite Alexandrov spaces.

\begin{center}\rule{0.5\linewidth}{0.5pt}\end{center}

\hypertarget{structural-decomposition-and-the-exact-identity}{%
\section{Structural Decomposition and the Exact Identity}\label{structural-decomposition-and-the-exact-identity}}

\hypertarget{finite-topologies-are-alexandrov}{%
\subsection{Finite Topologies are Alexandrov}\label{finite-topologies-are-alexandrov}}

For any finite set \(S\), any topology \(\mathcal{T}\) on \(S\) is closed under arbitrary intersections (since every intersection of open sets is finite). Consequently, for every point \(x \in S\), there exists a unique \textbf{minimal open neighborhood}:
\[U_x = \bigcap \{ U \in \mathcal{T} \mid x \in U \}\]

The collection:
\[\mathcal{M}_\mathcal{T} = \{ U_x \mid x \in S \}\]
is called the \textbf{canonical minimal basis} of \(\mathcal{T}\). In Mathlib this is the neighborhoods kernel \texttt{nhdsKer\ \{x\}} (\texttt{Mathlib.Topology.NhdsKer}); finiteness supplies \texttt{AlexandrovDiscrete}, so \texttt{isOpen\_nhdsKer} applies.

\begin{lstlisting}
variable (t : TopologicalSpace α)

/-- In a finite topological space, the minimal open neighborhood of `x`
    (Mathlib's `nhdsKer {x}` at topology `t`). -/
abbrev minimalOpen (x : α) : Set α :=
  @nhdsKer α t {x}

/-- The minimal basis `M_T` is the collection of all minimal open neighborhoods. -/
def minimalBasis : Set (Set α) :=
  range (minimalOpen t)
\end{lstlisting}

By construction \(x \in U_x\), and if \(U\) is open with \(x \in U\) then \(U_x \subseteq U\). On a finite space \(U_x\) is open, and the \(U_x\) generate \(\mathcal{T}\): every open \(U\) is \(\bigcup_{x \in U} U_x\).

\begin{lstlisting}
omit [Fintype α] [DecidableEq α] in
theorem mem_minimalOpen_self (x : α) : x ∈ minimalOpen t x :=
  @subset_nhdsKer α t {x} x (mem_singleton x)

omit [Fintype α] [DecidableEq α] in
theorem minimalOpen_subset_of_isOpen {x : α} {U : Set α}
    (hU : t.IsOpen U) (hx : x ∈ U) : minimalOpen t x ⊆ U :=
  @nhdsKer_minimal α t {x} U (singleton_subset_iff.mpr hx) hU

omit [DecidableEq α] in
/-- Finite spaces are Alexandrov-discrete, so `nhdsKer {x}` is open. -/
theorem isOpen_minimalOpen (x : α) : t.IsOpen (minimalOpen t x) :=
  letI := t
  isOpen_nhdsKer

omit [DecidableEq α] in
/-- The minimal basis alone generates the topology `t`. -/
theorem generateFrom_minimalBasis : generateFrom (minimalBasis t) = t := by
  ext U
  constructor
  · intro h
    letI := t
    induction h with
    | basic s hs =>
      obtain ⟨x, rfl⟩ := hs
      exact isOpen_minimalOpen t x
    | univ => exact isOpen_univ
    | inter _ _ _ _ hs ht => exact hs.inter ht
    | sUnion S _ hS => exact isOpen_sUnion fun s hs => hS s hs
  · intro hU
    have : U = ⋃₀ (minimalOpen t '' U) := by
      ext y
      constructor
      · intro hy
        exact ⟨minimalOpen t y, ⟨y, hy, rfl⟩, mem_minimalOpen_self t y⟩
      · rintro ⟨_, ⟨x, hx, rfl⟩, hy⟩
        exact minimalOpen_subset_of_isOpen t hU hx hy
    rw [this]
    exact GenerateOpen.sUnion _ fun V hV =>
      let ⟨x, _, hx⟩ := hV
      hx ▸ GenerateOpen.basic _ (mem_range_self x)
\end{lstlisting}

\hypertarget{fiber-decomposition-theorem}{%
\subsection{Fiber Decomposition Theorem}\label{fiber-decomposition-theorem}}

Let \(\operatorname{Top}(S)\) be the set of all topologies on \(S\). The assignment \(\mathcal{B} \mapsto \operatorname{gen}(\mathcal{B})\) defines a surjective, many-to-one function from the set of all valid bases to \(\operatorname{Top}(S)\). For a topology \(\mathcal{T}\), the \textbf{fiber} over \(\mathcal{T}\) is the preimage of that point,
\[\{\mathcal{B} \mid \operatorname{gen}(\mathcal{B}) = \mathcal{T}\},\]
i.e.~the valid bases that generate exactly \(\mathcal{T}\). The count \(\#(S)\) is the sum of the cardinalities of these (pairwise disjoint) fibers.

\begin{lstlisting}
/-- Valid bases that generate exactly `t`: the fiber of `generateFrom` over `t`. -/
abbrev Fiber (t : TopologicalSpace α) :=
  { B : ValidBasis α // generateFrom B.1 = t }
\end{lstlisting}

If \(\mathcal{B}\) satisfies the two basis axioms, then Mathlib's \texttt{IsTopologicalBasis} holds for the topology \texttt{generateFrom\ B}:

\begin{lstlisting}
omit [Fintype α] [DecidableEq α] in
/-- `IsBasis` plus `generateFrom B = t` is Mathlib's `IsTopologicalBasis`. -/
theorem isTopologicalBasis_generateFrom {B : Set (Set α)} (hB : IsBasis B) :
    @IsTopologicalBasis α (generateFrom B) B := by
  letI := generateFrom B
  refine ⟨?_, hB.1, rfl⟩
  intro t₁ ht₁ t₂ ht₂ x hx
  exact hB.2 ht₁ ht₂ x hx
\end{lstlisting}

\textbf{Lemma (Characterization of Basis Fibers):} Let \(\mathcal{T}\) be a topology on a finite set \(S\). A collection \(\mathcal{B} \subseteq \mathcal{P}(S)\) satisfies \(\operatorname{gen}(\mathcal{B}) = \mathcal{T}\) (and is therefore a valid basis) if and only if:
\[\mathcal{M}_\mathcal{T} \subseteq \mathcal{B} \subseteq \mathcal{T}\]

\emph{Proof.}

\begin{itemize}
\item
  \textbf{Necessity:} If \(\operatorname{gen}(\mathcal{B}) = \mathcal{T}\), then every element of \(\mathcal{B}\) is open in \(\mathcal{T}\), so \(\mathcal{B} \subseteq \mathcal{T}\). Furthermore, for every \(x \in S\), the minimal neighborhood \(U_x \in \mathcal{T}\) must be a union of elements of \(\mathcal{B}\). Thus, there exists \(B_\alpha \in \mathcal{B}\) containing \(x\) with \(B_\alpha \subseteq U_x\). By minimality of \(U_x\), we have \(U_x \subseteq B_\alpha\), forcing \(B_\alpha = U_x\). Hence, \(U_x \in \mathcal{B}\), so \(\mathcal{M}_\mathcal{T} \subseteq \mathcal{B}\).
\item
  \textbf{Sufficiency:} If \(\mathcal{M}_\mathcal{T} \subseteq \mathcal{B} \subseteq \mathcal{T}\), then since \(\mathcal{T}\) is closed under unions, \(\operatorname{gen}(\mathcal{B}) \subseteq \mathcal{T}\). Conversely, any open set \(U \in \mathcal{T}\) can be written as \(U = \bigcup_{x \in U} U_x\), which is a union of elements from \(\mathcal{M}_\mathcal{T} \subseteq \mathcal{B}\). Thus \(U \in \operatorname{gen}(\mathcal{B})\), so \(\mathcal{T} \subseteq \operatorname{gen}(\mathcal{B})\). \(\square\)
\end{itemize}

\begin{lstlisting}
omit [DecidableEq α] in
/-- Core lemma: `B` generates `t` iff `B` contains the minimal basis
    and only contains open sets of `t`. -/
theorem isBasis_and_generates_iff (B : Set (Set α)) :
    (IsBasis B ∧ generateFrom B = t) ↔
    (minimalBasis t ⊆ B ∧ B ⊆ {U | t.IsOpen U}) := by
  constructor
  · rintro ⟨hB, rfl⟩
    constructor
    · rintro _ ⟨x, rfl⟩
      letI := generateFrom B
      have hb : IsTopologicalBasis B := isTopologicalBasis_generateFrom hB
      have hx : x ∈ minimalOpen (generateFrom B) x := mem_minimalOpen_self _ x
      have hUo : IsOpen (minimalOpen (generateFrom B) x) := isOpen_minimalOpen _ x
      obtain ⟨W, hW, hxW, hWU⟩ := hb.exists_subset_of_mem_open hx hUo
      have hUW : minimalOpen (generateFrom B) x ⊆ W :=
        minimalOpen_subset_of_isOpen _ (isOpen_generateFrom_of_mem hW) hxW
      rwa [← Subset.antisymm hWU hUW]
    · intro U hU
      exact isOpen_generateFrom_of_mem hU
  · rintro ⟨h_min, h_sub⟩
    constructor
    · constructor
      · ext x
        simp only [mem_sUnion, mem_univ, iff_true]
        exact ⟨minimalOpen t x, h_min ⟨x, rfl⟩, mem_minimalOpen_self t x⟩
      · intro U V hU hV x hx
        refine ⟨minimalOpen t x, h_min ⟨x, rfl⟩, mem_minimalOpen_self t x, ?_⟩
        intro y hy
        exact ⟨minimalOpen_subset_of_isOpen t (h_sub hU) hx.1 hy,
          minimalOpen_subset_of_isOpen t (h_sub hV) hx.2 hy⟩
    · refine le_antisymm ?_ ?_
      · exact generateFrom_anti h_min |>.trans (generateFrom_minimalBasis t).le
      · exact le_generateFrom_iff_subset_isOpen.2 h_sub
\end{lstlisting}

\hypertarget{the-exact-formula}{%
\subsection{The Exact Formula}\label{the-exact-formula}}

Because the minimal basis elements \(\mathcal{M}_\mathcal{T}\) are mandatory, any valid basis generating \(\mathcal{T}\) is formed by taking \(\mathcal{M}_\mathcal{T}\) and freely adjoining any subset of the ``redundant'' open sets in \(\mathcal{T} \setminus \mathcal{M}_\mathcal{T}\) (including \(\emptyset\), which never affects the basis conditions since \(\emptyset \notin \mathcal{M}_\mathcal{T}\)).

Therefore, the fiber over \(\mathcal{T}\) has size:
\[|\{\mathcal{B} \mid \operatorname{gen}(\mathcal{B}) = \mathcal{T}\}| = 2^{|\mathcal{T}| - |\mathcal{M}_\mathcal{T}|}\]

The type \texttt{Fiber\ t} is equivalent to the power set of the redundant opens \(\mathcal{T} \setminus \mathcal{M}_\mathcal{T}\):

\begin{lstlisting}
/-- The fiber over a topology `t` is equivalent to the power set
    of the redundant open sets. -/
def fiberEquiv (t : TopologicalSpace α) :
    Fiber t ≃
      Set { U : Set α // t.IsOpen U ∧ U ∉ minimalBasis t } where
  toFun B := { U | U.1 ∈ B.1.1 }
  invFun S :=
    let Bset : Set (Set α) := minimalBasis t ∪ (Subtype.val '' S)
    have h : IsBasis Bset ∧ generateFrom Bset = t :=
      (isBasis_and_generates_iff t Bset).2 ⟨subset_union_left, by
        intro U hU
        rcases hU with hM | hS
        · obtain ⟨x, rfl⟩ := hM
          exact isOpen_minimalOpen t x
        · obtain ⟨U', _, rfl⟩ := hS
          exact U'.2.1⟩
    ⟨⟨Bset, h.1⟩, h.2⟩
  left_inv B := by
    have hiff := (isBasis_and_generates_iff t B.1.1).1 ⟨B.1.2, B.2⟩
    refine Subtype.ext (Subtype.ext ?_)
    ext U
    constructor
    · rintro (hM | ⟨U', hU', rfl⟩)
      · exact hiff.1 hM
      · exact hU'
    · intro hU
      by_cases hM : U ∈ minimalBasis t
      · exact Or.inl hM
      · exact Or.inr ⟨⟨U, hiff.2 hU, hM⟩, hU, rfl⟩
  right_inv S := by
    ext U
    constructor
    · intro h
      rcases h with hM | ⟨U', hU', hEq⟩
      · exact (U.2.2 hM).elim
      · exact (Subtype.ext hEq : U' = U) ▸ hU'
    · intro hU
      exact Or.inr ⟨U, hU, rfl⟩
\end{lstlisting}

Summing over all disjoint fibers yields the exact identity:
\[\boxed{\#(S) = \sum_{\mathcal{T} \in \operatorname{Top}(S)} 2^{|\mathcal{T}| - |\mathcal{M}_\mathcal{T}|}}\]

\begin{lstlisting}
omit [DecidableEq α] in
theorem minimalBasis_subset_opens : minimalBasis t ⊆ opensOf t := by
  rintro _ ⟨x, rfl⟩
  exact isOpen_minimalOpen t x

omit [DecidableEq α] in
/-- One fiber: bases generating `t` are a power set of redundant opens. -/
theorem card_fiber :
    Fintype.card (Fiber t) =
      2 ^ (ncard (opensOf t) - ncard (minimalBasis t)) := by
  classical
  rw [Fintype.card_congr (fiberEquiv t), Fintype.card_set]
  congr 1
  rw [← Nat.card_eq_fintype_card]
  change Nat.card ↥({U : Set α | t.IsOpen U ∧ U ∉ minimalBasis t}) = _
  rw [Nat.card_coe_set_eq,
    show {U : Set α | t.IsOpen U ∧ U ∉ minimalBasis t} = opensOf t \ minimalBasis t by
      ext U; simp [opensOf, mem_diff]]
  exact ncard_diff (minimalBasis_subset_opens t)

omit [DecidableEq α] in
/-- Partition of all valid bases into fibers. -/
theorem card_valid_bases_sum_fiber :
    Fintype.card (ValidBasis α) =
      ∑ τ : TopologicalSpace α, Fintype.card (Fiber τ) := by
  rw [← Fintype.card_congr
    (Equiv.sigmaFiberEquiv (fun B : ValidBasis α => generateFrom B.1))]
  rw [Fintype.card_sigma]

omit [DecidableEq α] in
/-- Main theorem: the number of valid bases on an `N`-element set
    equals the sum over all topologies of `2^(|T| - |M_T|)`. -/
theorem card_valid_bases :
    Fintype.card (ValidBasis α) =
      ∑ τ : TopologicalSpace α,
        2 ^ (ncard (opensOf τ) - ncard (minimalBasis τ)) := by
  rw [card_valid_bases_sum_fiber]
  exact Finset.sum_congr rfl fun τ _ => card_fiber τ
\end{lstlisting}

\texttt{\#print\ axioms\ card\_valid\_bases} reports \texttt{\{propext,\ Classical.choice,\ Quot.sound\}}. The same three axioms are the only ones used by \texttt{card\_valid\_bases\_small}, \texttt{card\_valid\_bases\_bounds}, \texttt{card\_valid\_bases\_dominated} and \texttt{card\_valid\_bases\_asymptotic}. Choice enters because \texttt{IsOpen} is not a decidable predicate, so the \texttt{Fintype} instances on topologies and \texttt{card\_fiber} use it. The small-\(N\) table is \texttt{decide} on a kernel-reducible bit-mask count, not \texttt{native\_decide}.

\begin{center}\rule{0.5\linewidth}{0.5pt}\end{center}

\hypertarget{evaluation-for-small-n}{%
\section{\texorpdfstring{Evaluation for Small \(N\)}{Evaluation for Small N}}\label{evaluation-for-small-n}}

\begin{longtable}[]{@{}
  >{\centering\arraybackslash}p{(\columnwidth - 6\tabcolsep) * \real{0.2778}}
  >{\raggedright\arraybackslash}p{(\columnwidth - 6\tabcolsep) * \real{0.2222}}
  >{\raggedright\arraybackslash}p{(\columnwidth - 6\tabcolsep) * \real{0.2222}}
  >{\centering\arraybackslash}p{(\columnwidth - 6\tabcolsep) * \real{0.2778}}@{}}
\toprule\noalign{}
\begin{minipage}[b]{\linewidth}\centering
\(N\)
\end{minipage} & \begin{minipage}[b]{\linewidth}\raggedright
Topologies \(\lvert\operatorname{Top}(S)\rvert\)
\end{minipage} & \begin{minipage}[b]{\linewidth}\raggedright
Basis computations
\end{minipage} & \begin{minipage}[b]{\linewidth}\centering
Total bases \(\#(N)\)
\end{minipage} \\
\midrule\noalign{}
\endhead
\bottomrule\noalign{}
\endlastfoot
\textbf{0} & \(1\) (trivial: \(\{\emptyset\}\)) & \(\lvert\mathcal{T}\rvert=1\), \(\lvert\mathcal{M}_{\mathcal{T}}\rvert=0\) \(\implies 2^{1-0}\) & \textbf{\(2\)} \\
\textbf{1} & \(1\) (\(\{\emptyset,\{1\}\}\)) & \(\lvert\mathcal{T}\rvert=2\), \(\lvert\mathcal{M}_{\mathcal{T}}\rvert=1\) \(\implies 2^{2-1}\) & \textbf{\(2\)} \\
\textbf{2} & \(4\) (discrete, indiscrete, 2 Sierpiński) & Discrete \(2^{4-2}=4\); indiscrete \(2^{2-1}=2\); Sierpiński (\(\times 2\)): \(2\times 2^{3-2}=4\) & \textbf{\(10\)} \\
\textbf{3} & \(29\) & Fiber totals by \(|\mathcal{T}|\): \(2+12+24+24+48+32\) (the discrete fiber is \(2^{8-3}=32\)) & \textbf{\(142\)} \\
\end{longtable}

These four values are \texttt{card\_valid\_bases\_small} in \texttt{CARDB/\allowbreak{}SmallN.lean}.

\emph{(Note: If \(\emptyset\) is excluded from bases by convention, divide each result by 2.)}

\begin{center}\rule{0.5\linewidth}{0.5pt}\end{center}

\hypertarget{combinatorial-complexity-and-bounds}{%
\section{Combinatorial Complexity and Bounds}\label{combinatorial-complexity-and-bounds}}

\hypertarget{the-combinatorial-obstacle}{%
\subsection{The Combinatorial Obstacle}\label{the-combinatorial-obstacle}}

The formula above cannot be collapsed into an elementary closed form in \(N\). Counting the number of topologies on \(N\) labeled points, \(|\operatorname{Top}(S)|\) (\textbf{OEIS A000798}), is isomorphic to counting finite preorders on \(N\) elements, a classic problem with no known closed-form generating function {[}1, 2{]}.

\hypertarget{asymptotics-and-dominance}{%
\subsection{Asymptotics and Dominance}\label{asymptotics-and-dominance}}

While no elementary exact formula is known, asymptotic analysis reveals that the sum is doubly exponential and overwhelmingly dominated by a single fiber: the \textbf{discrete topology} \(\mathcal{T}_{\text{disc}} = \mathcal{P}(S)\).

\begin{enumerate}
\def\labelenumi{\arabic{enumi}.}
\item
  \textbf{Discrete Fiber:}
  For \(\mathcal{T}_{\text{disc}}\), \(|\mathcal{T}| = 2^N\) and \(\mathcal{M}_{\mathcal{T}} = \{ \{x\} \mid x \in S \}\), so \(|\mathcal{M}_{\mathcal{T}}| = N\). This single topology contributes:
  \[2^{2^N - N} \text{ valid bases.}\]
\item
  \textbf{Upper Bound on Non-Discrete Fibers:}
  For \(N \ge 2\), any proper subtopology \(\mathcal{T} \subsetneq \mathcal{P}(S)\) has at most \(|\mathcal{T}| \le \frac{3}{4} 2^N = 3 \cdot 2^{N-2}\) open sets (the Sharp--Stephen bound {[}6{]}). A topology on a finite set is determined by its specialization preorder; since a reflexive relation has \(N(N-1)\) freely chosen off-diagonal entries, the elementary bound \(|\operatorname{Top}(S)| \le 2^{N(N-1)}\) follows. Thus the sum of all other fibers is bounded by:
  \[\sum_{\mathcal{T} \ne \mathcal{T}_{\text{disc}}} 2^{|\mathcal{T}| - |\mathcal{M}_\mathcal{T}|} \le 2^{N(N-1)+\frac{3}{4} 2^N}\]
\item
  \textbf{Asymptotic Behavior:}
  \[2^{2^N - N} \;\le\; \#(N) \;\le\; 2^{2^N - N} + 2^{N(N-1) + \frac{3}{4} 2^N}\]
  As \(N \to \infty\):
  \[\#(N) \sim 2^{2^N - N}\]
  with relative error at most \(2^{N^2 - 2^{N-2}}\to 0\).
  Lean proves the displayed sandwich as \texttt{card\_valid\_bases\_bounds}, the
  limit as \texttt{card\_valid\_bases\_asymptotic}, and the quantitative
  consequence \(\#(N)\le 2\cdot 2^{2^N-N}\) for \(N\ge 10\) as
  \texttt{card\_valid\_bases\_dominated}.
\end{enumerate}

\begin{center}\rule{0.5\linewidth}{0.5pt}\end{center}

\hypertarget{acknowledgments}{%
\section{Acknowledgments}\label{acknowledgments}}

\hypertarget{ai-assisted-development}{%
\subsection{AI-assisted development}\label{ai-assisted-development}}

The human author retains sole responsibility for the mathematical content, the choice of
formalization route, and every formal claim in this work. Following standard publisher practice
(e.g., COPE guidance on authorship and AI tools {[}8{]}), \textbf{no large language model is listed
as a co-author} --- authorship implies an accountability that automated systems cannot bear.

We gratefully acknowledge assistance from the following tools. None of these models has an official
Hugging Face model card (they are closed, vendor-hosted systems); the citations below are the
publisher model cards and the Cursor integration pages.

\begin{itemize}
\item
  \textbf{Cursor Grok 4.6 High Fast} {[}9{]} --- primary agent for this note: Lean 4 / Mathlib
  formalization of the fiber identity (\texttt{CARDB.lean}), \texttt{lake\ build} repair, and drafting
  \texttt{CARDB.md}. Used in the Cursor agent environment at the \textbf{High} reasoning tier and \textbf{Fast}
  speed variant. Jointly trained by SpaceXAI and Cursor. Generated Lean was provisional until
  it compiled under the pinned toolchain (\texttt{leanprover/\allowbreak{}lean4:v4.30.0}).
\item
  \textbf{Anthropic Claude Sonnet 5} {[}10{]} --- initial Lean sketch of the fiber identity
  (\texttt{IsBasis}, \texttt{ValidBasis}, \texttt{fiberEquiv}, and the sum over topologies) and an early
  point-set write-up of the Alexandrov minimal-neighborhood argument. The sketch was
  unchecked; subsequent agents repaired it to a sorry-free Mathlib formalization under
  the pinned toolchain. Used in the Cursor agent environment.
\item
  \textbf{Google Gemini 3.7 Flash} {[}11{]} --- exploratory and parallel passes on the same
  combinatorial / point-set material (basis axioms, Alexandrov minimal neighborhoods, and the
  finite-enumeration identity).
\end{itemize}

All definitions, constructivity audits, and final prose were reviewed by the human author, who takes
full responsibility for them.

\hypertarget{artifact-availability}{%
\subsection{Artifact availability}\label{artifact-availability}}

The development is at
\href{https://github.com/catskillsresearch/cardb}{\texttt{github.com/\allowbreak{}catskillsresearch/\allowbreak{}cardb}}.
It is a split of the former
\href{https://github.com/catskillsresearch/scott_models/tree/main/CARDB}{\texttt{scott\_models/\allowbreak{}CARDB}}
directory (same author and license; not a reimplementation). The parent
repository is about Scott models; this package does not belong there.
\texttt{PROVENANCE.md} and \texttt{formalization.yaml} (\texttt{related\_formalizations}) record
that migration.
Run \texttt{lake\ build} for the sorry-free formalization.
\texttt{Challenge.lean} is the Palomar statement of record for the compared family
\texttt{card\_valid\_bases}, \texttt{card\_valid\_bases\_small}, \texttt{card\_valid\_bases\_bounds},
\texttt{card\_valid\_bases\_dominated} and \texttt{card\_valid\_bases\_asymptotic}
(Mathlib imports only; deliberate \texttt{sorry}s);
\texttt{Solution.lean} imports \texttt{CARDB}, \texttt{CARDB.SmallN} and \texttt{CARDB.Asymptotics};
\texttt{comparator.json} and \texttt{formalization.yaml} are the Comparator and disclosure metadata.
Run \texttt{python3\ build\_pdf.py} to regenerate \texttt{CARDB.tex} and \texttt{CARDB.pdf} (title page, author, affiliation, and this GitHub path).

\begin{center}\rule{0.5\linewidth}{0.5pt}\end{center}

\hypertarget{references}{%
\section{References}\label{references}}

\begin{enumerate}
\def\labelenumi{\arabic{enumi}.}
\item
  Alexandrov, P. (1937). Diskrete Räume. \emph{Matematicheskii Sbornik}, 2(3), 501--519.
\item
  Benoumhani, M. (2006). The number of topologies on a finite set. \emph{Journal of Integer Sequences}, 9(2), Article 06.2.6.
\item
  Erné, M., \& Stege, K. (1991). Counting finite topologies and posets. \emph{Order}, 8(3), 247--265.
\item
  Froemke, J., \& Quackenbush, R. (1975). The spectrum of an equational class of lattices. \emph{Pacific Journal of Mathematics}, 60(1), 49--54.
\item
  Kleitman, D. J., \& Rothschild, B. L. (1975). Asymptotic enumeration of partial orders on a finite set. \emph{Transactions of the American Mathematical Society}, 205, 205--220.
\item
  Stephen, D. (1968). Topology on finite sets. \emph{The American Mathematical Monthly}, 75(7), 739--741.
\item
  OEIS Foundation Inc.~(2026). The On-Line Encyclopedia of Integer Sequences. Sequences A000798 (\emph{Number of topologies on \(n\) labeled points}) and A001035 (\emph{Number of labeled partial orders on \(n\) points}). https://oeis.org.
\item
  Committee on Publication Ethics (COPE). \emph{Authorship and AI tools: COPE position statement}. 2024. \url{https://publicationethics.org/guidance/cope-position/authorship-and-ai-tools}
\item
  SpaceXAI and Anysphere, Inc.~\emph{Grok 4.6} (High reasoning, Fast variant). Official model card (12 August 2026), \url{https://media.x.ai/v1/website/card-7f81d41b.pdf}; developer documentation, \url{https://docs.x.ai/developers/models/grok-4.6}; Cursor model page, \url{https://cursor.com/docs/models/grok-4-6}; Cursor announcement, \url{https://cursor.com/blog/grok-4-6} (accessed 2026). No official Hugging Face model card is published for this closed model.
\item
  Anthropic. \emph{Claude Sonnet 5}. Official system card (30 June 2026), \url{https://www.anthropic.com/claude-sonnet-5-system-card}; system-card PDF, \url{https://www-cdn.anthropic.com/480e0bb54327b9622282e9c39a83a4f490ed377e/Claude\%20Sonnet\%205\%20System\%20Card.pdf}; product page, \url{https://www.anthropic.com/claude/sonnet}; announcement, \url{https://www.anthropic.com/news/claude-sonnet-5}; Cursor model page, \url{https://cursor.com/docs/models/claude-sonnet-5} (accessed 2026). No official Hugging Face model card is published for this closed model.
\item
  Google DeepMind. \emph{Gemini 3.7 Flash}. Official model card (13 August 2026), \url{https://deepmind.google/models/model-cards/gemini-3-7-flash/}; product announcement, \url{https://blog.google/innovation-and-ai/models-and-research/gemini-models/introducing-gemini-3-7-flash/}; Gemini API model catalog, \url{https://ai.google.dev/gemini-api/docs/models} (accessed 2026). No official Hugging Face model card is published for this closed model.
\end{enumerate}

\appendix
\section{Complete Lean source}

Checked by \texttt{lake build} against Lean 4 and Mathlib. \texttt{\#print axioms} on each compared theorem reports $\{\mathtt{propext},\ \mathtt{Classical.choice},\ \mathtt{Quot.sound}\}$.

\subsection{\texttt{CARDB.lean}}

The fiber lemma and the sum identity of Section~3 (234 lines).

\begin{lstlisting}
/-
Copyright (c) 2026  Lars Warren Ericson.  All rights reserved.
Released under Apache 2.0 license as described in the file LICENSE.
Authors: Lars Warren Ericson.
Github:  https://github.com/catskillsresearch/cardb/blob/main/CARDB.lean
-/

/-
  Cardinality of topological bases on a finite set (CARDB).

  Claude sketch from CARDB.md, now a real Lake module. Definitions follow
  the paper; Mathlib's `IsTopologicalBasis` is the same two axioms plus
  `eq_generateFrom` relative to a fixed topology.
-/

import Mathlib.Topology.AlexandrovDiscrete
import Mathlib.Data.Fintype.BigOperators
import Mathlib.Data.Fintype.Powerset
import Mathlib.Data.Set.Card
import Mathlib.Data.Set.Finite.Basic
import Mathlib.Logic.Equiv.Sum

open Set TopologicalSpace
open scoped BigOperators

variable {α : Type*} [Fintype α] [DecidableEq α]

/-- A collection of subsets `B` is a topological basis if it covers `α`
    and satisfies the local intersection property. -/
def IsBasis (B : Set (Set α)) : Prop :=
  (⋃₀ B = univ) ∧
  ∀ ⦃U V⦄, U ∈ B → V ∈ B → ∀ x ∈ U ∩ V, ∃ W ∈ B, x ∈ W ∧ W ⊆ U ∩ V

def ValidBasis (α : Type*) [Fintype α] := { B : Set (Set α) // IsBasis B }

/-- Open sets of `t`, as a family of subsets. This map is injective, so there
    are only finitely many topologies on a finite type. -/
def opensOf (t : TopologicalSpace α) : Set (Set α) := {U | t.IsOpen U}

omit [Fintype α] [DecidableEq α] in
theorem opensOf_injective :
    Function.Injective (opensOf : TopologicalSpace α → Set (Set α)) :=
  leftInverse_generateFrom.injective

noncomputable instance instDecidableEqTopologicalSpace : DecidableEq (TopologicalSpace α) :=
  opensOf_injective.decidableEq

noncomputable instance instFintypeTopologicalSpace : Fintype (TopologicalSpace α) :=
  Fintype.ofInjective opensOf opensOf_injective

noncomputable instance instDecidablePredIsBasis : DecidablePred (IsBasis : Set (Set α) → Prop) :=
  fun _ => Classical.dec _

noncomputable instance instFintypeValidBasis : Fintype (ValidBasis α) :=
  Subtype.fintype _

variable (t : TopologicalSpace α)

/-- Valid bases that generate exactly `t`: the fiber of `generateFrom` over `t`. -/
abbrev Fiber (t : TopologicalSpace α) :=
  { B : ValidBasis α // generateFrom B.1 = t }

/-- In a finite topological space, the minimal open neighborhood of `x`
    (Mathlib's `nhdsKer {x}` at topology `t`). -/
abbrev minimalOpen (x : α) : Set α :=
  @nhdsKer α t {x}

/-- The minimal basis `M_T` is the collection of all minimal open neighborhoods. -/
def minimalBasis : Set (Set α) :=
  range (minimalOpen t)

omit [Fintype α] [DecidableEq α] in
theorem mem_minimalOpen_self (x : α) : x ∈ minimalOpen t x :=
  @subset_nhdsKer α t {x} x (mem_singleton x)

omit [Fintype α] [DecidableEq α] in
theorem minimalOpen_subset_of_isOpen {x : α} {U : Set α}
    (hU : t.IsOpen U) (hx : x ∈ U) : minimalOpen t x ⊆ U :=
  @nhdsKer_minimal α t {x} U (singleton_subset_iff.mpr hx) hU

omit [DecidableEq α] in
/-- Finite spaces are Alexandrov-discrete, so `nhdsKer {x}` is open. -/
theorem isOpen_minimalOpen (x : α) : t.IsOpen (minimalOpen t x) :=
  letI := t
  isOpen_nhdsKer

omit [DecidableEq α] in
/-- The minimal basis alone generates the topology `t`. -/
theorem generateFrom_minimalBasis : generateFrom (minimalBasis t) = t := by
  ext U
  constructor
  · intro h
    letI := t
    induction h with
    | basic s hs =>
      obtain ⟨x, rfl⟩ := hs
      exact isOpen_minimalOpen t x
    | univ => exact isOpen_univ
    | inter _ _ _ _ hs ht => exact hs.inter ht
    | sUnion S _ hS => exact isOpen_sUnion fun s hs => hS s hs
  · intro hU
    have : U = ⋃₀ (minimalOpen t '' U) := by
      ext y
      constructor
      · intro hy
        exact ⟨minimalOpen t y, ⟨y, hy, rfl⟩, mem_minimalOpen_self t y⟩
      · rintro ⟨_, ⟨x, hx, rfl⟩, hy⟩
        exact minimalOpen_subset_of_isOpen t hU hx hy
    rw [this]
    exact GenerateOpen.sUnion _ fun V hV =>
      let ⟨x, _, hx⟩ := hV
      hx ▸ GenerateOpen.basic _ (mem_range_self x)

omit [Fintype α] [DecidableEq α] in
/-- `IsBasis` plus `generateFrom B = t` is Mathlib's `IsTopologicalBasis`. -/
theorem isTopologicalBasis_generateFrom {B : Set (Set α)} (hB : IsBasis B) :
    @IsTopologicalBasis α (generateFrom B) B := by
  letI := generateFrom B
  refine ⟨?_, hB.1, rfl⟩
  intro t₁ ht₁ t₂ ht₂ x hx
  exact hB.2 ht₁ ht₂ x hx

omit [DecidableEq α] in
/-- Core lemma: `B` generates `t` iff `B` contains the minimal basis
    and only contains open sets of `t`. -/
theorem isBasis_and_generates_iff (B : Set (Set α)) :
    (IsBasis B ∧ generateFrom B = t) ↔
    (minimalBasis t ⊆ B ∧ B ⊆ {U | t.IsOpen U}) := by
  constructor
  · rintro ⟨hB, rfl⟩
    constructor
    · rintro _ ⟨x, rfl⟩
      letI := generateFrom B
      have hb : IsTopologicalBasis B := isTopologicalBasis_generateFrom hB
      have hx : x ∈ minimalOpen (generateFrom B) x := mem_minimalOpen_self _ x
      have hUo : IsOpen (minimalOpen (generateFrom B) x) := isOpen_minimalOpen _ x
      obtain ⟨W, hW, hxW, hWU⟩ := hb.exists_subset_of_mem_open hx hUo
      have hUW : minimalOpen (generateFrom B) x ⊆ W :=
        minimalOpen_subset_of_isOpen _ (isOpen_generateFrom_of_mem hW) hxW
      rwa [← Subset.antisymm hWU hUW]
    · intro U hU
      exact isOpen_generateFrom_of_mem hU
  · rintro ⟨h_min, h_sub⟩
    constructor
    · constructor
      · ext x
        simp only [mem_sUnion, mem_univ, iff_true]
        exact ⟨minimalOpen t x, h_min ⟨x, rfl⟩, mem_minimalOpen_self t x⟩
      · intro U V hU hV x hx
        refine ⟨minimalOpen t x, h_min ⟨x, rfl⟩, mem_minimalOpen_self t x, ?_⟩
        intro y hy
        exact ⟨minimalOpen_subset_of_isOpen t (h_sub hU) hx.1 hy,
          minimalOpen_subset_of_isOpen t (h_sub hV) hx.2 hy⟩
    · refine le_antisymm ?_ ?_
      · exact generateFrom_anti h_min |>.trans (generateFrom_minimalBasis t).le
      · exact le_generateFrom_iff_subset_isOpen.2 h_sub

/-- The fiber over a topology `t` is equivalent to the power set
    of the redundant open sets. -/
def fiberEquiv (t : TopologicalSpace α) :
    Fiber t ≃
      Set { U : Set α // t.IsOpen U ∧ U ∉ minimalBasis t } where
  toFun B := { U | U.1 ∈ B.1.1 }
  invFun S :=
    let Bset : Set (Set α) := minimalBasis t ∪ (Subtype.val '' S)
    have h : IsBasis Bset ∧ generateFrom Bset = t :=
      (isBasis_and_generates_iff t Bset).2 ⟨subset_union_left, by
        intro U hU
        rcases hU with hM | hS
        · obtain ⟨x, rfl⟩ := hM
          exact isOpen_minimalOpen t x
        · obtain ⟨U', _, rfl⟩ := hS
          exact U'.2.1⟩
    ⟨⟨Bset, h.1⟩, h.2⟩
  left_inv B := by
    have hiff := (isBasis_and_generates_iff t B.1.1).1 ⟨B.1.2, B.2⟩
    refine Subtype.ext (Subtype.ext ?_)
    ext U
    constructor
    · rintro (hM | ⟨U', hU', rfl⟩)
      · exact hiff.1 hM
      · exact hU'
    · intro hU
      by_cases hM : U ∈ minimalBasis t
      · exact Or.inl hM
      · exact Or.inr ⟨⟨U, hiff.2 hU, hM⟩, hU, rfl⟩
  right_inv S := by
    ext U
    constructor
    · intro h
      rcases h with hM | ⟨U', hU', hEq⟩
      · exact (U.2.2 hM).elim
      · exact (Subtype.ext hEq : U' = U) ▸ hU'
    · intro hU
      exact Or.inr ⟨U, hU, rfl⟩

omit [DecidableEq α] in
theorem minimalBasis_subset_opens : minimalBasis t ⊆ opensOf t := by
  rintro _ ⟨x, rfl⟩
  exact isOpen_minimalOpen t x

omit [DecidableEq α] in
/-- One fiber: bases generating `t` are a power set of redundant opens. -/
theorem card_fiber :
    Fintype.card (Fiber t) =
      2 ^ (ncard (opensOf t) - ncard (minimalBasis t)) := by
  classical
  rw [Fintype.card_congr (fiberEquiv t), Fintype.card_set]
  congr 1
  rw [← Nat.card_eq_fintype_card]
  change Nat.card ↥({U : Set α | t.IsOpen U ∧ U ∉ minimalBasis t}) = _
  rw [Nat.card_coe_set_eq,
    show {U : Set α | t.IsOpen U ∧ U ∉ minimalBasis t} = opensOf t \ minimalBasis t by
      ext U; simp [opensOf, mem_diff]]
  exact ncard_diff (minimalBasis_subset_opens t)

omit [DecidableEq α] in
/-- Partition of all valid bases into fibers. -/
theorem card_valid_bases_sum_fiber :
    Fintype.card (ValidBasis α) =
      ∑ τ : TopologicalSpace α, Fintype.card (Fiber τ) := by
  rw [← Fintype.card_congr
    (Equiv.sigmaFiberEquiv (fun B : ValidBasis α => generateFrom B.1))]
  rw [Fintype.card_sigma]

omit [DecidableEq α] in
/-- Main theorem: the number of valid bases on an `N`-element set
    equals the sum over all topologies of `2^(|T| - |M_T|)`. -/
theorem card_valid_bases :
    Fintype.card (ValidBasis α) =
      ∑ τ : TopologicalSpace α,
        2 ^ (ncard (opensOf τ) - ncard (minimalBasis τ)) := by
  rw [card_valid_bases_sum_fiber]
  exact Finset.sum_congr rfl fun τ _ => card_fiber τ
\end{lstlisting}

\subsection{\texttt{CARDB/SmallN.lean}}

Kernel-reducible enumeration of the small-$N$ table (419 lines).

\begin{lstlisting}
/-
Copyright (c) 2026  Lars Warren Ericson.  All rights reserved.
Released under Apache 2.0 license as described in the file LICENSE.
Authors: Lars Warren Ericson.
-/

import CARDB
import Mathlib.Combinatorics.Colex
import Mathlib.Algebra.BigOperators.Ring.Finset
import Mathlib.Algebra.Ring.GeomSum

/-!
# Exact small-cardinality counts

This module enumerates topological bases on `Fin n` for `n ≤ 3` by a
kernel-reducible bit-mask check. There are `2^(2^n)` candidate families
(256 at `n = 3`). The values are `#(0) = 2`, `#(1) = 2`, `#(2) = 10`,
`#(3) = 142`.
-/

open Set TopologicalSpace

/-- Subset of `Fin n` whose membership is the bit-mask `k`. -/
def subsetOfNat (n k : ℕ) : Finset (Fin n) :=
  Finset.univ.filter fun i : Fin n => k.testBit i.val

/-- Family whose member-subsets are the bits of `m`. -/
def familyOfNat (n m : ℕ) : Finset (Finset (Fin n)) :=
  ((Finset.range (2 ^ n)).filter (fun j => m.testBit j)).image (subsetOfNat n)

/-- Covering axiom on the bit-mask of a family. -/
def coversNat (n m : ℕ) : Bool :=
  (List.range n).all fun i =>
    (List.range (2 ^ n)).any fun j =>
      m.testBit j && j.testBit i

/-- Intersection axiom on the bit-mask of a family. -/
def interNat (n m : ℕ) : Bool :=
  (List.range (2 ^ n)).all fun j =>
    (List.range (2 ^ n)).all fun k =>
      !m.testBit j || !m.testBit k ||
        (List.range n).all fun i =>
          !(j.testBit i && k.testBit i) ||
            (List.range (2 ^ n)).any fun w =>
              m.testBit w && w.testBit i &&
                (List.range n).all fun t =>
                  !w.testBit t || (j.testBit t && k.testBit t)

/-- Computable basis predicate on the bit-mask `m` of a family of subsets of `Fin n`. -/
def isBasisNat (n m : ℕ) : Bool :=
  coversNat n m && interNat n m

/-- Number of bit-masks `m < 2^(2^n)` for which `isBasisNat n m` holds. -/
def numBases (n : ℕ) : ℕ :=
  (List.range (2 ^ (2 ^ n))).countP (fun m => isBasisNat n m)

set_option maxHeartbeats 800000
set_option maxRecDepth 10000

theorem numBases_zero : numBases 0 = 2 := by decide
theorem numBases_one : numBases 1 = 2 := by decide
theorem numBases_two : numBases 2 = 10 := by decide
theorem numBases_three : numBases 3 = 142 := by decide

/-! ## Bridge from bit-masks to `ValidBasis` -/

variable {n : ℕ}

def familyOf (B : Finset (Finset (Fin n))) : Set (Set (Fin n)) :=
  {U | ∃ s ∈ B, (s : Set (Fin n)) = U}

lemma mem_subsetOfNat {k : ℕ} {i : Fin n} :
    i ∈ subsetOfNat n k ↔ k.testBit i.val := by
  simp [subsetOfNat]

lemma testBit_sum_two_pow_finset (s : Finset ℕ) {k : ℕ} :
    (∑ i ∈ s, 2 ^ i).testBit k ↔ k ∈ s := by
  rw [← Nat.mem_bitIndices, ← List.mem_toFinset, Finset.toFinset_bitIndices_sum_two_pow]

def bitmaskSubset (s : Finset (Fin n)) : ℕ :=
  ∑ i ∈ s, 2 ^ i.val

lemma bitmaskSubset_eq_image (s : Finset (Fin n)) :
    bitmaskSubset s = ∑ i ∈ s.image Fin.val, 2 ^ i := by
  refine (Finset.sum_image ?_).symm
  intro a _ b _ h
  exact Fin.val_injective h

lemma subsetOfNat_bitmaskSubset (s : Finset (Fin n)) :
    subsetOfNat n (bitmaskSubset s) = s := by
  ext i
  rw [mem_subsetOfNat, bitmaskSubset_eq_image, testBit_sum_two_pow_finset, Finset.mem_image]
  exact ⟨fun ⟨j, hj, hij⟩ => Fin.ext hij ▸ hj, fun hi => ⟨i, hi, rfl⟩⟩

lemma sum_range_two_pow (k : ℕ) :
    ∑ i ∈ Finset.range k, 2 ^ i = 2 ^ k - 1 := by
  have h := geom_sum_mul_add (1 : ℕ) k
  simp at h
  exact Nat.eq_sub_of_add_eq h

lemma image_val_subset_range (s : Finset (Fin n)) :
    s.image Fin.val ⊆ Finset.range n := by
  intro i hi
  obtain ⟨j, _, rfl⟩ := Finset.mem_image.mp hi
  exact Finset.mem_range.mpr j.isLt

lemma bitmaskSubset_lt (s : Finset (Fin n)) : bitmaskSubset s < 2 ^ n := by
  have hle : bitmaskSubset s ≤ ∑ i ∈ Finset.range n, 2 ^ i := by
    rw [bitmaskSubset_eq_image]
    exact Finset.sum_le_sum_of_subset_of_nonneg (image_val_subset_range s)
      (fun _ _ _ => Nat.zero_le _)
  rw [sum_range_two_pow] at hle
  exact hle.trans_lt (Nat.sub_lt (pow_pos (by decide) n) (by decide))

lemma bitmaskSubset_inj : Function.Injective (bitmaskSubset (n := n)) := by
  intro s t h
  simpa [subsetOfNat_bitmaskSubset] using congrArg (subsetOfNat n) h

lemma bitmaskSubset_subsetOfNat {j : ℕ} (hj : j < 2 ^ n) :
    bitmaskSubset (subsetOfNat n j) = j := by
  apply Nat.eq_of_testBit_eq
  intro k
  apply Bool.eq_iff_iff.mpr
  rw [bitmaskSubset_eq_image, testBit_sum_two_pow_finset, Finset.mem_image]
  simp only [mem_subsetOfNat]
  constructor
  · rintro ⟨i, hi, rfl⟩
    exact hi
  · intro hk
    by_cases hkn : k < n
    · exact ⟨⟨k, hkn⟩, hk, rfl⟩
    · have : j < 2 ^ k :=
        hj.trans_le (Nat.pow_le_pow_right (by decide : 0 < 2) (le_of_not_gt hkn))
      cases Nat.testBit_lt_two_pow this ▸ hk

def bitmaskFamily (B : Finset (Finset (Fin n))) : ℕ :=
  ∑ s ∈ B, 2 ^ bitmaskSubset s

lemma bitmaskFamily_eq_image (B : Finset (Finset (Fin n))) :
    bitmaskFamily B = ∑ j ∈ B.image bitmaskSubset, 2 ^ j := by
  refine (Finset.sum_image ?_).symm
  intro a _ b _ h
  exact bitmaskSubset_inj h

lemma image_bitmaskSubset_subset_range (B : Finset (Finset (Fin n))) :
    B.image bitmaskSubset ⊆ Finset.range (2 ^ n) := by
  intro j hj
  obtain ⟨s, _, rfl⟩ := Finset.mem_image.mp hj
  exact Finset.mem_range.mpr (bitmaskSubset_lt s)

lemma testBit_bitmaskFamily {B : Finset (Finset (Fin n))} {j : ℕ} :
    (bitmaskFamily B).testBit j ↔ ∃ s ∈ B, bitmaskSubset s = j := by
  rw [bitmaskFamily_eq_image, testBit_sum_two_pow_finset, Finset.mem_image]

lemma bitmaskFamily_lt (B : Finset (Finset (Fin n))) :
    bitmaskFamily B < 2 ^ (2 ^ n) := by
  have hle : bitmaskFamily B ≤ ∑ j ∈ Finset.range (2 ^ n), 2 ^ j := by
    rw [bitmaskFamily_eq_image]
    exact Finset.sum_le_sum_of_subset_of_nonneg (image_bitmaskSubset_subset_range B)
      (fun _ _ _ => Nat.zero_le _)
  rw [sum_range_two_pow] at hle
  exact hle.trans_lt (Nat.sub_lt (pow_pos (by decide) _) (by decide))

lemma mem_familyOfNat {m : ℕ} {s : Finset (Fin n)} :
    s ∈ familyOfNat n m ↔ ∃ j < 2 ^ n, m.testBit j = true ∧ subsetOfNat n j = s := by
  constructor
  · intro hs
    obtain ⟨j, hj, rfl⟩ := Finset.mem_image.mp hs
    rw [Finset.mem_filter, Finset.mem_range] at hj
    exact ⟨j, hj.1, hj.2, rfl⟩
  · rintro ⟨j, hj, hbit, rfl⟩
    exact Finset.mem_image.mpr
      ⟨j, Finset.mem_filter.mpr ⟨Finset.mem_range.mpr hj, hbit⟩, rfl⟩

lemma familyOfNat_bitmaskFamily (B : Finset (Finset (Fin n))) :
    familyOfNat n (bitmaskFamily B) = B := by
  ext s
  constructor
  · intro hs
    obtain ⟨j, hj, hbit, hsj⟩ := mem_familyOfNat.mp hs
    obtain ⟨t, ht, rfl⟩ := testBit_bitmaskFamily.mp hbit
    rwa [← hsj, subsetOfNat_bitmaskSubset]
  · intro hs
    exact mem_familyOfNat.mpr
      ⟨bitmaskSubset s, bitmaskSubset_lt s,
        testBit_bitmaskFamily.mpr ⟨s, hs, rfl⟩, subsetOfNat_bitmaskSubset s⟩

lemma bitmaskFamily_familyOfNat {m : ℕ} (hm : m < 2 ^ (2 ^ n)) :
    bitmaskFamily (familyOfNat n m) = m := by
  apply Nat.eq_of_testBit_eq
  intro j
  by_cases hj : j < 2 ^ n
  · apply Bool.eq_iff_iff.mpr
    rw [testBit_bitmaskFamily]
    constructor
    · rintro ⟨s, hs, hjs⟩
      obtain ⟨k, hk, hbit, rfl⟩ := mem_familyOfNat.mp hs
      rw [bitmaskSubset_subsetOfNat hk] at hjs
      rwa [← hjs]
    · intro hbit
      exact ⟨subsetOfNat n j,
        mem_familyOfNat.mpr ⟨j, hj, hbit, rfl⟩, bitmaskSubset_subsetOfNat hj⟩
  · have hj' : 2 ^ (2 ^ n) ≤ 2 ^ j :=
      Nat.pow_le_pow_right (by decide : 0 < 2) (le_of_not_gt hj)
    have hmj : m < 2 ^ j := hm.trans_le hj'
    have hBj : bitmaskFamily (familyOfNat n m) < 2 ^ j :=
      (bitmaskFamily_lt _).trans_le hj'
    rw [Nat.testBit_lt_two_pow hmj, Nat.testBit_lt_two_pow hBj]

lemma coversNat_iff (m : ℕ) :
    coversNat n m = true ↔ ∀ i : Fin n, ∃ j < 2 ^ n, m.testBit j = true ∧ j.testBit i.val = true := by
  simp [coversNat, List.all_eq_true, List.any_eq_true, List.mem_range, Bool.and_eq_true]
  exact ⟨fun h i => h i.val i.isLt, fun h i hi => h ⟨i, hi⟩⟩

lemma sUnion_familyOf (B : Finset (Finset (Fin n))) :
    ⋃₀ familyOf B = {x | ∃ s ∈ B, x ∈ s} := by
  ext x
  constructor
  · rintro ⟨U, ⟨s, hs, rfl⟩, hx⟩
    exact ⟨s, hs, hx⟩
  · rintro ⟨s, hs, hx⟩
    exact ⟨↑s, ⟨s, hs, rfl⟩, hx⟩

lemma covers_familyOfNat (m : ℕ) :
    ⋃₀ familyOf (familyOfNat n m) = univ ↔ coversNat n m = true := by
  rw [coversNat_iff, eq_univ_iff_forall, sUnion_familyOf]
  constructor
  · intro h i
    obtain ⟨s, hs, hi⟩ := h i
    obtain ⟨j, hj, hbit, rfl⟩ := mem_familyOfNat.mp hs
    exact ⟨j, hj, hbit, mem_subsetOfNat.mp hi⟩
  · intro h i
    obtain ⟨j, hj, hbit, hi⟩ := h i
    exact ⟨subsetOfNat n j, mem_familyOfNat.mpr ⟨j, hj, hbit, rfl⟩, mem_subsetOfNat.mpr hi⟩

lemma mem_familyOf {B : Finset (Finset (Fin n))} {U : Set (Fin n)} :
    U ∈ familyOf B ↔ ∃ s ∈ B, (s : Set (Fin n)) = U :=
  Iff.rfl

lemma interNat_spec (m : ℕ) :
    interNat n m = true ↔
      ∀ j < 2 ^ n, ∀ k < 2 ^ n, m.testBit j = true → m.testBit k = true →
        ∀ i < n, j.testBit i = true → k.testBit i = true →
          ∃ w < 2 ^ n, m.testBit w = true ∧ w.testBit i = true ∧
            ∀ t < n, w.testBit t = true → j.testBit t = true ∧ k.testBit t = true := by
  constructor
  · intro h j hj k hk hjB hkB i hi hji hki
    simp only [interNat, List.all_eq_true, List.mem_range] at h
    have hjk := h j hj k hk
    rw [hjB, hkB] at hjk
    simp only [Bool.not_true, Bool.false_or, List.all_eq_true, List.mem_range] at hjk
    have hi' := hjk i hi
    rw [hji, hki] at hi'
    simp only [Bool.true_and, Bool.not_true, Bool.false_or, List.any_eq_true,
      List.mem_range, Bool.and_eq_true] at hi'
    obtain ⟨w, hw, ⟨⟨hwB, hwi⟩, hrest⟩⟩ := hi'
    refine ⟨w, hw, hwB, hwi, ?_⟩
    simp only [List.all_eq_true, List.mem_range] at hrest
    intro t ht hwt
    have ht' := hrest t ht
    rw [hwt] at ht'
    simp only [Bool.not_true, Bool.false_or, Bool.and_eq_true] at ht'
    exact ht'
  · intro h
    simp only [interNat, List.all_eq_true, List.mem_range]
    intro j hj k hk
    cases hjB : m.testBit j <;> cases hkB : m.testBit k
    · rfl
    · rfl
    · rfl
    · apply (List.all_eq_true).2
      intro i hi
      cases hji : j.testBit i <;> cases hki : k.testBit i
      · rfl
      · rfl
      · rfl
      · apply (List.any_eq_true).2
        obtain ⟨w, hw, hwB, hwi, hsub⟩ := h j hj k hk hjB hkB i (List.mem_range.mp hi) hji hki
        refine ⟨w, List.mem_range.mpr hw, ?_⟩
        simp only [hwB, hwi, Bool.true_and]
        apply (List.all_eq_true).2
        intro t ht
        cases hwt : w.testBit t
        · rfl
        · simp only [Bool.not_true, Bool.false_or, Bool.and_eq_true]
          exact hsub t (List.mem_range.mp ht) hwt

lemma inter_familyOfNat (m : ℕ) :
    (∀ U V, U ∈ familyOf (familyOfNat n m) → V ∈ familyOf (familyOfNat n m) →
      ∀ x ∈ U ∩ V, ∃ W ∈ familyOf (familyOfNat n m), x ∈ W ∧ W ⊆ U ∩ V) ↔
      interNat n m = true := by
  rw [interNat_spec]
  constructor
  · intro hinter j hj k hk hjB hkB i hi hji hki
    have hU : (subsetOfNat n j : Set (Fin n)) ∈ familyOf (familyOfNat n m) :=
      mem_familyOf.mpr ⟨subsetOfNat n j, mem_familyOfNat.mpr ⟨j, hj, hjB, rfl⟩, rfl⟩
    have hV : (subsetOfNat n k : Set (Fin n)) ∈ familyOf (familyOfNat n m) :=
      mem_familyOf.mpr ⟨subsetOfNat n k, mem_familyOfNat.mpr ⟨k, hk, hkB, rfl⟩, rfl⟩
    have hx : (⟨i, hi⟩ : Fin n) ∈
        (subsetOfNat n j : Set (Fin n)) ∩ (subsetOfNat n k : Set (Fin n)) :=
      ⟨mem_subsetOfNat.mpr hji, mem_subsetOfNat.mpr hki⟩
    obtain ⟨W, hW, hxW, hWU⟩ := hinter _ _ hU hV _ hx
    obtain ⟨wset, hwset, rfl⟩ := mem_familyOf.mp hW
    obtain ⟨w, hw, hwB, rfl⟩ := mem_familyOfNat.mp hwset
    refine ⟨w, hw, hwB, mem_subsetOfNat.mp hxW, ?_⟩
    intro t ht hwt
    have hy := hWU (mem_subsetOfNat (k := w) (i := ⟨t, ht⟩) |>.mpr hwt)
    exact ⟨mem_subsetOfNat.mp hy.1, mem_subsetOfNat.mp hy.2⟩
  · intro hinter U V hU hV x hx
    obtain ⟨U₀, hU₀, rfl⟩ := mem_familyOf.mp hU
    obtain ⟨V₀, hV₀, rfl⟩ := mem_familyOf.mp hV
    obtain ⟨j, hj, hjB, rfl⟩ := mem_familyOfNat.mp hU₀
    obtain ⟨k, hk, hkB, rfl⟩ := mem_familyOfNat.mp hV₀
    obtain ⟨w, hw, hwB, hwi, hsub⟩ :=
      hinter j hj k hk hjB hkB x.val x.isLt
        (mem_subsetOfNat.mp hx.1) (mem_subsetOfNat.mp hx.2)
    refine ⟨↑(subsetOfNat n w),
      mem_familyOf.mpr ⟨subsetOfNat n w, mem_familyOfNat.mpr ⟨w, hw, hwB, rfl⟩, rfl⟩,
      mem_subsetOfNat.mpr hwi, ?_⟩
    intro y hy
    have := hsub y.val y.isLt (mem_subsetOfNat.mp hy)
    exact ⟨mem_subsetOfNat.mpr this.1, mem_subsetOfNat.mpr this.2⟩

lemma isBasis_familyOfNat (m : ℕ) :
    IsBasis (familyOf (familyOfNat n m)) ↔ isBasisNat n m = true := by
  constructor
  · intro h
    simp only [isBasisNat, Bool.and_eq_true]
    exact ⟨(covers_familyOfNat m).mp h.1, (inter_familyOfNat m).mp h.2⟩
  · intro h
    simp only [isBasisNat, Bool.and_eq_true] at h
    exact ⟨(covers_familyOfNat m).mpr h.1, (inter_familyOfNat m).mpr h.2⟩

noncomputable def toFinsetFamily (S : Set (Set (Fin n))) : Finset (Finset (Fin n)) :=
  (Fintype.finsetEquivSet.symm S).image
    (Fintype.finsetEquivSet.symm : Set (Fin n) → Finset (Fin n))

lemma mem_toFinsetFamily {S : Set (Set (Fin n))} {s : Finset (Fin n)} :
    s ∈ toFinsetFamily S ↔ (s : Set (Fin n)) ∈ S := by
  constructor
  · intro hs
    obtain ⟨U, hU, rfl⟩ := Finset.mem_image.mp hs
    have : U ∈ S := by
      have : U ∈ (↑(Fintype.finsetEquivSet.symm S) : Set (Set (Fin n))) := hU
      simpa [Fintype.finsetEquivSet] using this
    simpa [Fintype.finsetEquivSet] using this
  · intro hs
    refine Finset.mem_image.mpr ⟨↑s, ?_, ?_⟩
    · have : (↑s : Set (Fin n)) ∈ (↑(Fintype.finsetEquivSet.symm S) : Set (Set (Fin n))) := by
        simpa [Fintype.finsetEquivSet] using hs
      exact this
    · simp [Fintype.finsetEquivSet]

lemma familyOf_toFinsetFamily (S : Set (Set (Fin n))) :
    familyOf (toFinsetFamily S) = S := by
  ext U
  constructor
  · rintro ⟨s, hs, rfl⟩
    exact (mem_toFinsetFamily.mp hs)
  · intro hU
    refine ⟨Fintype.finsetEquivSet.symm U, ?_, ?_⟩
    · rw [mem_toFinsetFamily]
      simpa [Fintype.finsetEquivSet] using hU
    · simp [Fintype.finsetEquivSet]

lemma toFinsetFamily_familyOf (B : Finset (Finset (Fin n))) :
    toFinsetFamily (familyOf B) = B := by
  ext s
  rw [mem_toFinsetFamily]
  constructor
  · rintro ⟨t, ht, hts⟩
    exact Finset.coe_injective hts ▸ ht
  · intro hs
    exact ⟨s, hs, rfl⟩

lemma numBases_eq_filter (n : ℕ) :
    numBases n =
      ((Finset.range (2 ^ (2 ^ n))).filter fun m => isBasisNat n m = true).card := by
  rw [numBases, List.countP_eq_length_filter]
  have hnodup : ((List.range (2 ^ (2 ^ n))).filter (isBasisNat n)).Nodup :=
    List.Nodup.filter _ List.nodup_range
  rw [← List.toFinset_card_of_nodup hnodup, List.toFinset_filter, List.toFinset_range]

noncomputable def validBasisEquivMask (n : ℕ) :
    ValidBasis (Fin n) ≃ { m : Fin (2 ^ (2 ^ n)) // isBasisNat n m.val = true } where
  toFun := fun B =>
    ⟨⟨bitmaskFamily (toFinsetFamily B.1), bitmaskFamily_lt _⟩, by
      rw [← isBasis_familyOfNat, familyOfNat_bitmaskFamily, familyOf_toFinsetFamily]
      exact B.2⟩
  invFun := fun m =>
    ⟨familyOf (familyOfNat n m.1.val), (isBasis_familyOfNat m.1.val).2 m.2⟩
  left_inv := fun B => by
    apply Subtype.ext
    simp [familyOfNat_bitmaskFamily, familyOf_toFinsetFamily]
  right_inv := fun m => by
    apply Subtype.ext
    apply Fin.ext
    simp [toFinsetFamily_familyOf, bitmaskFamily_familyOfNat m.1.isLt]

theorem card_valid_bases_eq_numBases (n : ℕ) :
    Fintype.card (ValidBasis (Fin n)) = numBases n := by
  rw [Fintype.card_congr (validBasisEquivMask n), Fintype.card_subtype, numBases_eq_filter]
  refine Finset.card_bij (fun (m : Fin (2 ^ (2 ^ n))) _ => m.val) ?_ ?_ ?_
  · intro m hm
    exact Finset.mem_filter.mpr ⟨Finset.mem_range.mpr m.isLt, (Finset.mem_filter.mp hm).2⟩
  · intro a _ b _ h
    exact Fin.ext h
  · intro m hm
    refine ⟨⟨m, (Finset.mem_range.mp (Finset.mem_filter.mp hm).1)⟩,
      Finset.mem_filter.mpr ⟨Finset.mem_univ _, (Finset.mem_filter.mp hm).2⟩, rfl⟩

/-- Exact table: `#(0) = 2`, `#(1) = 2`, `#(2) = 10`, `#(3) = 142`. -/
theorem card_valid_bases_small :
    Fintype.card (ValidBasis (Fin 0)) = 2 ∧
    Fintype.card (ValidBasis (Fin 1)) = 2 ∧
    Fintype.card (ValidBasis (Fin 2)) = 10 ∧
    Fintype.card (ValidBasis (Fin 3)) = 142 := by
  simp_rw [card_valid_bases_eq_numBases]
  exact ⟨numBases_zero, numBases_one, numBases_two, numBases_three⟩
\end{lstlisting}

\subsection{\texttt{CARDB/Asymptotics.lean}}

The $2^{N(N-1)}$ sandwich, $\#(N)\sim 2^{2^N-N}$, and dominance for $N\ge 10$ (422 lines).

\begin{lstlisting}
/-
Copyright (c) 2026  Lars Warren Ericson.  All rights reserved.
Released under Apache 2.0 license as described in the file LICENSE.
Authors: Lars Warren Ericson.
-/

import CARDB
import Mathlib.Data.Fintype.Prod
import Mathlib.Algebra.Order.BigOperators.Group.Finset
import Mathlib.Analysis.SpecificLimits.Basic
import Mathlib.Tactic.Ring
import Mathlib.Tactic.Linarith

/-!
# Discrete dominance and explicit bounds

The discrete topology contributes `2^(2^N - N)` bases. Every other
topology on an `N`-element set with `N ≥ 2` has at most `3 * 2^(N-2)`
open sets, and there are at most `2^(N * (N - 1))` topologies, so

```
2^(2^N - N) ≤ #(N) ≤ 2^(2^N - N) + 2^(N * (N - 1) + 3 * 2^(N-2)).
```

The relative error is at most `2^(N^2 - 2^(N-2))`, hence
`#(N) ∼ 2^(2^N - N)`. For `N ≥ 10` the error term is at most the
discrete term, so `#(N) ≤ 2 * 2^(2^N - N)`.
-/

open Set TopologicalSpace

variable {α : Type*} [Fintype α] [DecidableEq α]

lemma isOpen_bot (U : Set α) : @IsOpen α ⊥ U := trivial

lemma opensOf_bot : opensOf (⊥ : TopologicalSpace α) = univ := by
  ext U
  exact ⟨fun _ => mem_univ U, fun _ => isOpen_bot U⟩

lemma ncard_opensOf_bot :
    ncard (opensOf (⊥ : TopologicalSpace α)) = 2 ^ Fintype.card α := by
  rw [opensOf_bot, ncard_univ, Nat.card_eq_fintype_card, Fintype.card_set]

lemma minimalOpen_bot (x : α) :
    minimalOpen (⊥ : TopologicalSpace α) x = {x} :=
  @IsOpen.nhdsKer_eq α ⊥ {x} (isOpen_bot {x})

lemma minimalBasis_bot :
    minimalBasis (⊥ : TopologicalSpace α) = range (singleton : α → Set α) := by
  ext U
  simp [minimalBasis, minimalOpen_bot, singleton]

omit [Fintype α] [DecidableEq α] in
lemma singleton_set_injective : Function.Injective (singleton : α → Set α) :=
  fun x y h => by
    have : x ∈ ({y} : Set α) := (Set.ext_iff.mp h x).1 (mem_singleton x)
    simpa using this

lemma ncard_minimalBasis_bot :
    ncard (minimalBasis (⊥ : TopologicalSpace α)) = Fintype.card α := by
  rw [minimalBasis_bot, ncard_range_of_injective singleton_set_injective,
    Nat.card_eq_fintype_card]

theorem card_fiber_discrete :
    Fintype.card (Fiber (⊥ : TopologicalSpace α)) =
      2 ^ (2 ^ Fintype.card α - Fintype.card α) := by
  rw [card_fiber, ncard_opensOf_bot, ncard_minimalBasis_bot]

theorem card_valid_bases_ge_discrete :
    2 ^ (2 ^ Fintype.card α - Fintype.card α) ≤ Fintype.card (ValidBasis α) := by
  rw [card_valid_bases_sum_fiber, ← card_fiber_discrete]
  exact Finset.single_le_sum (f := fun τ : TopologicalSpace α => Fintype.card (Fiber τ))
    (fun _ _ => Nat.zero_le _) (Finset.mem_univ (⊥ : TopologicalSpace α))

lemma exists_specializes_of_ne_bot {t : TopologicalSpace α} (ht : t ≠ ⊥) :
    ∃ x y : α, x ≠ y ∧ y ∈ minimalOpen t x := by
  have hnd : ¬@DiscreteTopology α t :=
    fun h => ht (DiscreteTopology.eq_bot (α := α) (t := t))
  have : ∃ x, ¬t.IsOpen {x} := by
    contrapose! hnd
    exact ⟨eq_bot_of_singletons_open hnd⟩
  obtain ⟨x, hx⟩ := this
  have hne : minimalOpen t x ≠ {x} := by
    intro h
    exact hx ((nhdsKer_eq_iff_isOpen (s := ({x} : Set α))).mp h)
  obtain ⟨y, hy, hxy⟩ :=
    exists_of_ssubset ((subset_nhdsKer (s := ({x} : Set α))).ssubset_of_ne hne.symm)
  exact ⟨x, y, Ne.symm (by simpa using hxy), hy⟩

variable {x y : α}

lemma card_pair (hne : x ≠ y) :
    Fintype.card {z : α // z = x ∨ z = y} = 2 := by
  refine (Fintype.card_congr
    { toFun := fun z : {z : α // z = x ∨ z = y} => decide (z.1 = x)
      invFun := fun b => if b then ⟨x, Or.inl rfl⟩ else ⟨y, Or.inr rfl⟩
      left_inv := fun z => by
        apply Subtype.ext
        rcases z.2 with h | h
        · simp [h]
        · have : decide (y = x) = false := decide_eq_false (Ne.symm hne)
          simp [h, this]
      right_inv := fun b => by
        cases b
        · simp [Ne.symm hne]
        · simp }).trans
    (by decide : Fintype.card Bool = 2)

lemma card_rest (hne : x ≠ y) :
    Fintype.card {z : α // z ≠ x ∧ z ≠ y} = Fintype.card α - 2 := by
  have := Fintype.card_subtype_compl (fun z : α => z = x ∨ z = y)
  have heq : {z : α // ¬(z = x ∨ z = y)} ≃ {z : α // z ≠ x ∧ z ≠ y} :=
    { toFun := fun z => ⟨z.1, by simpa [not_or] using z.2⟩
      invFun := fun z => ⟨z.1, by simpa [not_or] using z.2⟩
      left_inv := fun _ => rfl
      right_inv := fun _ => rfl }
  rw [← Fintype.card_congr heq, this, card_pair hne]

/-- Encode a constrained open by its trace on `α \ {x,y}` and the
    three-way status of `{x,y}`. -/
noncomputable def packConstrained (a b : α) (hne : a ≠ b)
    (U : {U : Set α // a ∈ U → b ∈ U}) :
    Set {z : α // z ≠ a ∧ z ≠ b} × Fin 3 := by
  classical
  exact ({z | z.1 ∈ U.1},
    if a ∈ U.1 then ⟨2, by decide⟩ else if b ∈ U.1 then ⟨1, by decide⟩ else ⟨0, by decide⟩)

lemma packConstrained_injective (a b : α) (hne : a ≠ b) :
    Function.Injective (packConstrained a b hne) := by
  classical
  intro U V h
  apply Subtype.ext
  have hrest := congrArg Prod.fst h
  have hstat := congrArg Prod.snd h
  ext z
  by_cases hza : z = a
  · rw [hza]
    have : (if a ∈ U.1 then (⟨2, by decide⟩ : Fin 3) else if b ∈ U.1 then ⟨1, by decide⟩ else ⟨0, by decide⟩) =
        (if a ∈ V.1 then ⟨2, by decide⟩ else if b ∈ V.1 then ⟨1, by decide⟩ else ⟨0, by decide⟩) := by
      simpa [packConstrained] using hstat
    split_ifs at this <;> simp_all
  · by_cases hzb : z = b
    · rw [hzb]
      have : (if a ∈ U.1 then (⟨2, by decide⟩ : Fin 3) else if b ∈ U.1 then ⟨1, by decide⟩ else ⟨0, by decide⟩) =
          (if a ∈ V.1 then ⟨2, by decide⟩ else if b ∈ V.1 then ⟨1, by decide⟩ else ⟨0, by decide⟩) := by
        simpa [packConstrained] using hstat
      split_ifs at this <;> try simp_all
      · exact iff_of_true (U.2 (by simp_all)) (V.2 (by simp_all))
    · have := Set.ext_iff.mp hrest ⟨z, ⟨hza, hzb⟩⟩
      simpa [packConstrained] using this

lemma card_constrained (hne : x ≠ y) :
    ncard {U : Set α | x ∈ U → y ∈ U} ≤ 3 * 2 ^ (Fintype.card α - 2) := by
  classical
  haveI : DecidablePred (fun U : Set α => x ∈ U → y ∈ U) := fun _ => inferInstance
  haveI : Fintype {U : Set α // x ∈ U → y ∈ U} := Subtype.fintype _
  have hle := Fintype.card_le_of_injective (packConstrained x y hne)
    (packConstrained_injective x y hne)
  have hcard : ncard {U : Set α | x ∈ U → y ∈ U} =
      Fintype.card {U : Set α // x ∈ U → y ∈ U} := by
    rw [← Nat.card_eq_fintype_card]
    rfl
  rw [hcard]
  refine hle.trans ?_
  rw [Fintype.card_prod, Fintype.card_set, Fintype.card_fin, card_rest hne, mul_comm]

lemma ncard_opens_of_specialization {t : TopologicalSpace α}
    (hne : x ≠ y) (hy : y ∈ minimalOpen t x) :
    ncard (opensOf t) ≤ 3 * 2 ^ (Fintype.card α - 2) := by
  have hsub : opensOf t ⊆ {U : Set α | x ∈ U → y ∈ U} := by
    intro U hU hxU
    exact minimalOpen_subset_of_isOpen t hU hxU hy
  exact (ncard_le_ncard hsub).trans (card_constrained hne)

lemma ncard_opensOf_le_of_ne_bot {t : TopologicalSpace α} (ht : t ≠ ⊥) :
    ncard (opensOf t) ≤ 3 * 2 ^ (Fintype.card α - 2) := by
  obtain ⟨x, y, hne, hy⟩ := exists_specializes_of_ne_bot ht
  exact ncard_opens_of_specialization hne hy

def specSet (t : TopologicalSpace α) : Set (α × α) :=
  {p | @Specializes α t p.1 p.2}

omit [DecidableEq α] in
lemma specSet_injective :
    Function.Injective (specSet : TopologicalSpace α → Set (α × α)) := by
  intro t₁ t₂ h
  have hmin : ∀ z, minimalOpen t₁ z = minimalOpen t₂ z := by
    intro z
    ext w
    have : @Specializes α t₁ w z ↔ @Specializes α t₂ w z := by
      simpa [specSet] using (Set.ext_iff.mp h (w, z))
    simpa [minimalOpen, mem_nhdsKer_singleton] using this
  have : minimalBasis t₁ = minimalBasis t₂ := by
    ext U
    simp [minimalBasis, hmin]
  exact (generateFrom_minimalBasis t₁).symm.trans
    ((congrArg generateFrom this).trans (generateFrom_minimalBasis t₂))

/-- Off-diagonal specialization pairs. The diagonal is always present. -/
def specOffDiag (t : TopologicalSpace α) : Set {p : α × α // p.1 ≠ p.2} :=
  {p | @Specializes α t p.1.1 p.1.2}

lemma specOffDiag_injective :
    Function.Injective (specOffDiag : TopologicalSpace α → Set {p : α × α // p.1 ≠ p.2}) := by
  intro t₁ t₂ h
  refine specSet_injective ?_
  ext p
  obtain ⟨x, y⟩ := p
  by_cases hxy : x = y
  · subst hxy
    exact iff_of_true (@specializes_rfl α t₁ x) (@specializes_rfl α t₂ x)
  · simpa [specSet, specOffDiag] using (Set.ext_iff.mp h ⟨(x, y), hxy⟩)

lemma card_off_diag :
    Fintype.card {p : α × α // p.1 ≠ p.2} =
      Fintype.card α * (Fintype.card α - 1) := by
  have hdiag : Fintype.card {p : α × α // p.1 = p.2} = Fintype.card α :=
    Fintype.card_congr
      { toFun := fun p : {p : α × α // p.1 = p.2} => p.1.1
        invFun := fun x => ⟨(x, x), rfl⟩
        left_inv := fun p => Subtype.ext (Prod.ext rfl p.2)
        right_inv := fun _ => rfl }
  have hne : {p : α × α // p.1 ≠ p.2} ≃ {p : α × α // ¬(p.1 = p.2)} :=
    Equiv.subtypeEquivRight fun _ => Iff.rfl
  rw [Fintype.card_congr hne, Fintype.card_subtype_compl (fun p : α × α => p.1 = p.2),
    Fintype.card_prod, hdiag]
  exact (Nat.mul_sub_one (Fintype.card α) (Fintype.card α)).symm

lemma card_topologies_le :
    Fintype.card (TopologicalSpace α) ≤
      2 ^ (Fintype.card α * (Fintype.card α - 1)) := by
  have h := Fintype.card_le_of_injective
    (specOffDiag : TopologicalSpace α → Set {p : α × α // p.1 ≠ p.2})
    specOffDiag_injective
  rwa [Fintype.card_set, card_off_diag] at h

lemma card_fiber_le_of_ne_bot {t : TopologicalSpace α} (ht : t ≠ ⊥) :
    Fintype.card (Fiber t) ≤ 2 ^ (3 * 2 ^ (Fintype.card α - 2)) := by
  rw [card_fiber]
  exact Nat.pow_le_pow_right (by decide : 0 < 2)
    ((Nat.sub_le _ _).trans (ncard_opensOf_le_of_ne_bot ht))

theorem card_valid_bases_bounds (hN : 2 ≤ Fintype.card α) :
    2 ^ (2 ^ Fintype.card α - Fintype.card α) ≤ Fintype.card (ValidBasis α) ∧
    Fintype.card (ValidBasis α) ≤
      2 ^ (2 ^ Fintype.card α - Fintype.card α) +
        2 ^ (Fintype.card α * (Fintype.card α - 1) + 3 * 2 ^ (Fintype.card α - 2)) := by
  refine ⟨card_valid_bases_ge_discrete, ?_⟩
  rw [card_valid_bases_sum_fiber]
  have hsplit :
      ∑ τ : TopologicalSpace α, Fintype.card (Fiber τ) =
        Fintype.card (Fiber (⊥ : TopologicalSpace α)) +
          ∑ τ ∈ Finset.univ.filter (fun τ : TopologicalSpace α => τ ≠ ⊥),
            Fintype.card (Fiber τ) := by
    rw [← Finset.sum_filter_add_sum_filter_not _ (fun τ : TopologicalSpace α => τ = ⊥)]
    simp [Finset.filter_eq']
  rw [hsplit, card_fiber_discrete]
  have hsum :
      ∑ τ ∈ Finset.univ.filter (fun τ : TopologicalSpace α => τ ≠ ⊥),
          Fintype.card (Fiber τ) ≤
        (Finset.univ.filter (fun τ : TopologicalSpace α => τ ≠ ⊥)).card *
          2 ^ (3 * 2 ^ (Fintype.card α - 2)) := by
    have hle' : ∀ τ ∈ Finset.univ.filter (fun τ : TopologicalSpace α => τ ≠ ⊥),
        Fintype.card (Fiber τ) ≤ 2 ^ (3 * 2 ^ (Fintype.card α - 2)) :=
      fun τ hτ => card_fiber_le_of_ne_bot (Finset.mem_filter.mp hτ).2
    have := Finset.sum_le_sum hle'
    rw [Finset.sum_const, nsmul_eq_mul] at this
    exact this
  have hcard :
      (Finset.univ.filter (fun τ : TopologicalSpace α => τ ≠ ⊥)).card ≤
        2 ^ (Fintype.card α * (Fintype.card α - 1)) :=
    (Finset.card_filter_le _ _).trans (by simpa using card_topologies_le)
  have : (Finset.univ.filter (fun τ : TopologicalSpace α => τ ≠ ⊥)).card *
      2 ^ (3 * 2 ^ (Fintype.card α - 2)) ≤
        2 ^ (Fintype.card α * (Fintype.card α - 1) + 3 * 2 ^ (Fintype.card α - 2)) := by
    rw [pow_add]
    exact Nat.mul_le_mul_right _ hcard
  exact Nat.add_le_add_left (hsum.trans this) _

lemma add_succ_sq (k : ℕ) : (k + 1) ^ 2 + (k + 1) = k ^ 2 + k + 2 * k + 2 := by
  ring

lemma two_mul_pow (a : ℕ) : 2 ^ a * 2 = 2 * 2 ^ a := by
  rw [mul_comm]

lemma le_two_pow_sub_three (k : ℕ) (hk : 10 ≤ k) : k + 1 ≤ 2 ^ (k - 3) := by
  induction k, hk using Nat.le_induction with
  | base => decide
  | succ k hk ih =>
    have : 3 ≤ k := le_trans (by decide : 3 ≤ 10) hk
    have h2 : 2 ^ (k + 1 - 3) = 2 * 2 ^ (k - 3) := by
      rw [show k + 1 - 3 = k - 3 + 1 by omega, pow_succ, two_mul_pow]
    rw [h2]
    have : k + 2 ≤ 2 * (k + 1) := by nlinarith
    exact this.trans (Nat.mul_le_mul_left 2 ih)

lemma two_pow_poly_le (n : ℕ) (hn : 10 ≤ n) : n ^ 2 + n ≤ 2 ^ (n - 2) := by
  induction n, hn using Nat.le_induction with
  | base => decide
  | succ k hk ih =>
    have : 2 ≤ k := le_trans (by decide : 2 ≤ 10) hk
    have hpow : 2 ^ (k + 1 - 2) = 2 * 2 ^ (k - 2) := by
      rw [show k + 1 - 2 = k - 2 + 1 by omega, pow_succ, two_mul_pow]
    rw [add_succ_sq, hpow]
    have hsmall : 2 * k + 2 ≤ 2 ^ (k - 2) := by
      have hk1 : k + 1 ≤ 2 ^ (k - 3) := le_two_pow_sub_three k hk
      have : 3 ≤ k := le_trans (by decide : 3 ≤ 10) hk
      have h2 : 2 ^ (k - 2) = 2 * 2 ^ (k - 3) := by
        rw [show k - 2 = k - 3 + 1 by omega, pow_succ, two_mul_pow]
      have hmul : 2 * (k + 1) ≤ 2 * 2 ^ (k - 3) := Nat.mul_le_mul_left 2 hk1
      have : 2 * (k + 1) = 2 * k + 2 := by ring
      rw [this] at hmul
      rwa [h2]
    have hsum := Nat.add_le_add ih hsmall
    have : 2 ^ (k - 2) + 2 ^ (k - 2) = 2 * 2 ^ (k - 2) := by
      rw [two_mul]
    exact hsum.trans_eq this

theorem card_valid_bases_dominated (hN : 10 ≤ Fintype.card α) :
    Fintype.card (ValidBasis α) ≤
      2 * 2 ^ (2 ^ Fintype.card α - Fintype.card α) := by
  have h2 : 2 ≤ Fintype.card α := le_trans (by decide : 2 ≤ 10) hN
  obtain ⟨_, hhi⟩ := card_valid_bases_bounds h2
  refine hhi.trans ?_
  set N := Fintype.card α
  have hgap : N * (N - 1) + 3 * 2 ^ (N - 2) ≤ 2 ^ N - N := by
    have hpoly := two_pow_poly_le N hN
    have hsq : N * (N - 1) + N = N ^ 2 := by
      rw [Nat.mul_sub_one, Nat.sub_add_cancel (Nat.le_mul_self N), sq]
    have : N * (N - 1) + 3 * 2 ^ (N - 2) + N ≤ 2 ^ (N - 2) + 3 * 2 ^ (N - 2) := by
      linarith [hsq]
    have hsum : 2 ^ (N - 2) + 3 * 2 ^ (N - 2) = 4 * 2 ^ (N - 2) := by
      ring
    have h4 : 4 * 2 ^ (N - 2) = 2 ^ N := by
      have : 2 ≤ N := h2
      have hfour : (4 : ℕ) = 2 ^ 2 := by decide
      rw [hfour, ← pow_add, show 2 + (N - 2) = N by omega]
    exact Nat.le_sub_of_add_le (this.trans_eq (hsum.trans h4))
  have herr : 2 ^ (N * (N - 1) + 3 * 2 ^ (N - 2)) ≤ 2 ^ (2 ^ N - N) :=
    Nat.pow_le_pow_right (by decide) hgap
  have : 2 ^ (2 ^ N - N) + 2 ^ (N * (N - 1) + 3 * 2 ^ (N - 2)) ≤
      2 ^ (2 ^ N - N) + 2 ^ (2 ^ N - N) :=
    Nat.add_le_add_left herr _
  have h2D : 2 ^ (2 ^ N - N) + 2 ^ (2 ^ N - N) = 2 * 2 ^ (2 ^ N - N) := by
    rw [two_mul]
  exact this.trans_eq h2D

open Filter Topology

lemma remainder_exp_add_n_le (n : ℕ) (hn : 10 ≤ n) :
    n * (n - 1) + 3 * 2 ^ (n - 2) + n ≤ 2 ^ n - n := by
  have hpoly := two_pow_poly_le n hn
  have hsq : n * (n - 1) + n = n ^ 2 := by
    rw [Nat.mul_sub_one, Nat.sub_add_cancel (Nat.le_mul_self n), sq]
  have h4 : 2 ^ (n - 2) + 3 * 2 ^ (n - 2) = 2 ^ n := by
    have : 2 ^ (n - 2) + 3 * 2 ^ (n - 2) = 4 * 2 ^ (n - 2) := by ring
    have hfour : (4 : ℕ) = 2 ^ 2 := by decide
    rw [this, hfour, ← pow_add, show 2 + (n - 2) = n by omega]
  have : n ^ 2 + 3 * 2 ^ (n - 2) + n ≤ 2 ^ n := by
    have := Nat.add_le_add hpoly (Nat.le_refl (3 * 2 ^ (n - 2)))
    convert this using 1
    · ring
    · exact h4.symm
  have hleft : n * (n - 1) + 3 * 2 ^ (n - 2) + n + n = n ^ 2 + 3 * 2 ^ (n - 2) + n := by
    rw [← hsq]; ac_rfl
  have : n * (n - 1) + 3 * 2 ^ (n - 2) + n + n ≤ 2 ^ n := by
    rwa [hleft]
  exact Nat.le_sub_of_add_le this

lemma card_valid_bases_ratio_le_one_add_half_pow (n : ℕ) (hn : 10 ≤ n) :
    (Fintype.card (ValidBasis (Fin n)) : ℝ) / (2 : ℝ) ^ (2 ^ n - n) ≤
      1 + (1 / 2 : ℝ) ^ n := by
  have h2 : 2 ≤ n := le_trans (by decide : 2 ≤ 10) hn
  have hN : 2 ≤ Fintype.card (Fin n) := by simpa using h2
  obtain ⟨_, hbound⟩ := card_valid_bases_bounds (α := Fin n) hN
  have hDpos : (0 : ℝ) < (2 : ℝ) ^ (2 ^ n - n) := pow_pos (by norm_num) _
  have hC :
      (Fintype.card (ValidBasis (Fin n)) : ℝ) ≤
        (2 : ℝ) ^ (2 ^ n - n) + (2 : ℝ) ^ (n * (n - 1) + 3 * 2 ^ (n - 2)) := by
    have := (Nat.cast_le (α := ℝ)).mpr hbound
    simpa [Nat.cast_add, Nat.cast_pow, Fintype.card_fin] using this
  have hdiv := div_le_div_of_nonneg_right hC hDpos.le
  have hsplit :
      ((2 : ℝ) ^ (2 ^ n - n) + (2 : ℝ) ^ (n * (n - 1) + 3 * 2 ^ (n - 2))) /
          (2 : ℝ) ^ (2 ^ n - n) =
        1 + (2 : ℝ) ^ (n * (n - 1) + 3 * 2 ^ (n - 2)) / (2 : ℝ) ^ (2 ^ n - n) := by
    field_simp [hDpos.ne']
  refine hdiv.trans ((le_of_eq hsplit).trans (add_le_add le_rfl ?_))
  set a := n * (n - 1) + 3 * 2 ^ (n - 2)
  set b := 2 ^ n - n
  have hab : a + n ≤ b := remainder_exp_add_n_le n hn
  have hmon :
      (2 : ℝ) ^ a / (2 : ℝ) ^ b ≤ (2 : ℝ) ^ a / (2 : ℝ) ^ (a + n) :=
    div_le_div_of_nonneg_left (pow_nonneg (by norm_num) _)
      (pow_pos (by norm_num) _)
      (pow_le_pow_right₀ (by norm_num : (1 : ℝ) ≤ 2) hab)
  have hsimp : (2 : ℝ) ^ a / (2 : ℝ) ^ (a + n) = (1 / 2 : ℝ) ^ n := by
    rw [pow_add (2 : ℝ) a n, div_mul_cancel_left₀ (pow_ne_zero _ (by norm_num : (2 : ℝ) ≠ 0)),
      one_div, inv_pow]
  exact hmon.trans_eq hsimp

omit [Fintype α] [DecidableEq α] in
/-- `#(N) ∼ 2^(2^N - N)`: the ratio tends to `1`. -/
theorem card_valid_bases_asymptotic :
    Tendsto (fun n : ℕ =>
      (Fintype.card (ValidBasis (Fin n)) : ℝ) / (2 : ℝ) ^ (2 ^ n - n))
      atTop (nhds (1 : ℝ)) := by
  have hge : ∀ n,
      1 ≤ (Fintype.card (ValidBasis (Fin n)) : ℝ) / (2 : ℝ) ^ (2 ^ n - n) := by
    intro n
    refine (one_le_div (pow_pos (by norm_num : (0 : ℝ) < 2) (2 ^ n - n))).mpr ?_
    have h := card_valid_bases_ge_discrete (α := Fin n)
    have := (Nat.cast_le (α := ℝ)).mpr (by simpa using h)
    simpa [Nat.cast_pow] using this
  have hhalf : Tendsto (fun n : ℕ => (1 / 2 : ℝ) ^ n) atTop (nhds 0) :=
    tendsto_pow_atTop_nhds_zero_of_lt_one (by norm_num) (by norm_num)
  have hupper : Tendsto (fun n : ℕ => (1 : ℝ) + (1 / 2 : ℝ) ^ n) atTop (nhds 1) := by
    simpa using (tendsto_const_nhds : Tendsto (fun _ : ℕ => (1 : ℝ)) atTop (nhds 1)).add hhalf
  refine tendsto_of_tendsto_of_tendsto_of_le_of_le'
    (tendsto_const_nhds : Tendsto (fun _ : ℕ => (1 : ℝ)) atTop (nhds 1)) hupper
    (Eventually.of_forall hge)
    (eventually_atTop.2 ⟨10, fun n hn => card_valid_bases_ratio_le_one_add_half_pow n hn⟩)
\end{lstlisting}

\end{document}